\documentclass{amsart}
\usepackage[american]{babel}
\usepackage{amsopn}
\usepackage{amsmath,amssymb}
\usepackage{graphicx}
\usepackage{float}
\usepackage{fancyhdr}
\numberwithin{equation}{section}
\def\R{{\mathbb R}}
\def\E{{\mathbb E}}
\def\N{{\mathbb N}}

\def\spn{{\rm span}}
\def\F{{\mathcal F}}
\def\ch{{\rm cosh}}
\def\sech{{\rm sech}}
\def\sh{{\rm sinh}}
\def\p{{\partial}}
\def\xL{{\rm L}}
\def\xW{{\rm W}}
\def\xCn{{\rm C}}
\def\s{{\sigma}}
\def\xdif{{\rm d}}
\def\xLtwo{{{\rm L}^2}}
\def\xHone{{{\rm H}^1}}
\def\xHtwo{{{\rm H}^2}}
\def\L2{{\rm L}^2(\R^d)}
\def\xLn{{\rm L}}
\def\xLone{{\rm L}^1}
\def\xHn{{\rm H}}
\def\xtr{{\rm tr}}
\def\H1{{\xHone(\R^d)}}
\def\indic{{\rm {\large 1}\hspace{-2.3pt}{\large
l}}}
\def\xLinfty{{\rm L}^{\infty}}
\def\xim{{\rm im}}
\def\xker{{\rm ker}}
\theoremstyle{plain}
\newtheorem{thrm}{Theorem}[section]
\newtheorem{lmm}[thrm]{Lemma}
\newtheorem{crllr}[thrm]{Corollary}
\newtheorem{prpstn}[thrm]{Proposition}
\theoremstyle{definition}

\newtheorem{rmrk}[thrm]{Remark}
\title{Large deviations and support results for nonlinear
Schr\"odinger equations with additive noise and applications}
\author{Eric GAUTIER$^{1,2}$}

\begin{document}

\maketitle \pagestyle{myheadings} \markboth{{\sc E. Gautier}}{{\sc
Large deviations and support for nonlinear Schr\"odinger equations
with additive noise}}

\begin{abstract} Sample path large deviations
for the laws of the solutions of stochastic nonlinear
Schr\"odinger equations when the noise converges to zero are
presented. The noise is a complex additive gaussian noise. It is
white in time and colored space wise. The solutions may be global
or blow-up in finite time, the two cases are distinguished. The
results are stated in trajectory spaces endowed with projective
limit topologies. In this setting, the support of the law of the
solution is also characterized. As a consequence, results on the
law of the blow-up time and asymptotics when the noise converges
to zero are obtained. An application to the transmission of
solitary waves in fiber optics is also given.\vspace{0.3cm}

\end{abstract}

\footnotetext[1]{Laboratoire de Statistique, CREST-INSEE, URA
D2200, 3 avenue Pierre Larousse, 92240 Malakoff, France}

\footnotetext[2]{Laboratoire de Processus Stochastiques et
Statistique, IRMAR, UMR 6625, Universit\'e de Rennes 1, Campus de
Beaulieu, 35042 Rennes cedex, France;
\email{eric.gautier@bretagne.ens-cachan.fr}\vspace{0.2cm}

\noindent{\em Key Words:} } \thispagestyle{empty}
\section{Introduction}\label{s1}
In the present article, the stochastic nonlinear Schr\"odinger
(NLS) equation with a power law nonlinearity and an additive noise
is studied. The deterministic equation occurs as a basic model in
many areas of physics: hydrodynamics, plasma physics, nonlinear
optics, molecular biology. It describes the propagation of waves
in media with both nonlinear and dispersive responses. It is an
idealized model and does not take into account many aspects such
as inhomogeneities, high order terms, thermal fluctuations,
external forces which may be modeled as a random excitation (see
\cite{BCIR1,BCIR2,FKLT,konotop-vasquez}). Propagation in random
media may also be considered. The resulting re-scaled equation is
a random perturbation of the dynamical system of the following
form:
\begin{equation}
\label{e1} i \frac{\partial}{\partial t} \psi - (\Delta \psi
+\lambda|\psi|^{2\s}\psi)=\xi, \quad x \in \R^d, \quad t \geq 0,
\quad \lambda=\pm 1,
\end{equation}
where $\xi$ is a complex valued space-time white noise with

correlation function, following the denomination used in
\cite{FKLT},

$$\E\left[\xi(t_1,x_1)\bar{\xi}(t_2,x_2)\right]=D
\delta_{t_1-t_2}\otimes\delta_{x_1-x_2}$$ $D$ is the noise
intensity and $\delta$ denotes the Dirac mass. When $\lambda=1$
the nonlinearity is called focusing, otherwise it is defocusing.
\vspace{0.3cm}

With the notations of section \ref{s2}, the unbounded operator
$-i\Delta$ on $\L2$ with domain $\xHtwo(\R^d)$ is skew-adjoint.
Stone's theorem gives thus that it generates a unitary group
$(S(t)=e^{-it\Delta})_{t\in\R}$. The Fourier transform gives that
this group is also unitary on every Sobolev space based on $\L2$.
Consequently, there is no smoothing effect in the Sobolev spaces.

We are thus unable to treat the space-time white noise and will
consider a complex valued centered gaussian noise, white in time
and colored space wise.\vspace{0.3cm}

In the present article, the formalism of stochastic evolution
equations in Banach spaces as presented in \cite{DPZ} is adopted.
This point of view is preferred to the field and martingale
measure stochastic integral approach, see \cite{WAL}, in order to
use a particular property of the group, namely
hyper-contractivity. The Strichartz estimates, presented in the
next section, show that some integrability property is gained
through time integration and "convolution" with the group. In this
setting, the gaussian noise is defined as the time derivative in
the sense of distributions of a $Q$-Wiener process
$(W(t))_{t\in[0,+\infty)}$ on $\H1$. Here $Q$ is the covariance
operator of the law of the $\H1-$random variable $W(1)$, which is
a centered gaussian measure. With the It\^o notations, the
stochastic evolution equation is written
\begin{equation}
\label{e2} i \xdif u -(\Delta u + \lambda |u|^{2\s}u)\xdif t =
\xdif W.
\end{equation}
The initial datum $u_0$ is a function of $\H1$.  We will consider
solutions of NLS that are weak solutions in the sense used in the
analysis of partial differential equations or equivalently mild
solutions which satisfy
\begin{equation}
\label{e5}u(t)=S(t)u_0-i\lambda\int_0^tS(t-s)(|u(s)|^{2\sigma}u(s))ds-i\int_0^tS(t-s)dW(s).
\end{equation}

The well posedness of the Cauchy problem associated to (\ref{e1})
in the deterministic case depends on the size of $\s$. If
$\s<\frac{2}{d}$, the nonlinearity is subcritical and the Cauchy
problem is globally well posed in $\L2$ or $\H1$. If
$\s=\frac{2}{d}$, critical nonlinearity, or
$\frac{2}{d}<\s<\frac{2}{d-2}$ when $d\geq 3$ or simply
$\s>\frac{2}{d}$ otherwise, supercritical nonlinearity, the Cauchy
problem is locally well posed in $\H1$, see \cite{kato}. In this
latter case, if the nonlinearity is defocusing, the solution is
global. In the focusing case, when the nonlinearity is critical or
supercritical, some initial data yield global solutions while it
is known that other initial data yield solutions which blow up in
finite time, see \cite{CAZ,SS}.\vspace{0.3cm}

In \cite{dBD1}, the $\H1$ results have been generalized to the
stochastic case and existence and uniqueness results are obtained
for the stochastic equation under the same conditions on $\s$.
Continuity with respect to the initial data and the perturbation
is proved. It is shown that the proof of global existence for a
defocusing nonlinearity or for a focusing nonlinearity with a
subcritical exponent, could be adapted in the stochastic case even
if the momentum $$M(u)=\|u\|_{\L2}$$ and hamiltonian
$$H(u)=\frac{1}{2}\int_{\mathbb{R}^d}|\nabla u|^2
dx-\frac{\lambda}{2\sigma+2}\int_{\mathbb{R}^d}|u|^{2\sigma+2}dx$$
are no longer conserved. For a focusing nonlinearity and critical

or supercritical exponents, the solution may blow-up in finite
time. The blow-up time is denoted by $\tau(\omega)$. It satisfies
either $\lim_{t\rightarrow\tau(\omega)}\|u(t)\|_{\H1}=+\infty$ or
$\tau(\omega)=+\infty$, even if the solution is obtained by a
fixed point argument in a ball of a space of more regular
functions than $\xCn([0,T];\H1).$\vspace{0.3cm}

In this article, we are interested in the law of the paths of the
random solution. When the noise converges to zero, continuity with

respect to the perturbation gives that the law converges to the
Dirac mass on the deterministic solution. In the following, a
large deviation result is shown. It gives the rate of convergence
to zero, on an exponential scale, of the probability that paths
are in sets that do not contain the deterministic solution. A
general result is stated for the case where blow-up in finite time
is possible and a second one for the particular case where the
solutions are global. Also, the stronger is the topology, the
sharper are the estimates. We will therefore take advantage of the
variety of spaces that can be considered for the fixed point
argument, due to the integrability property, and present the large
deviation principles in trajectory spaces endowed with projective
limit topologies. A characterization of the support of the law of
the solution in these trajectory spaces is proved. The two results
can be transferred to weaker topologies or more generally by any
continuous mapping. The first application is a proof that, for
certain noises, with positive probability some solutions blow up
after any time $T$. Some estimates on the law of the blow-up time
when the noise converges to zero are also obtained. This study is
yet another contribution to the study of the influence of a noise

on the blow-up of the solutions of the focusing supercritical NLS,
see in the case of an additive noise \cite{dBD2,dBD4}. A second
application is given. It consists in obtaining similar results as
in \cite{FKLT} with an approach based on large deviations. The aim
is to compute estimates of error probability in signal
transmission in optical fibers when the medium is random and
nonlinear, for small noises.\vspace{0.3cm}

Section \ref{s2} is devoted to notations and properties of the
group, of the noise and of the stochastic convolution. An
extension of the result of continuity with respect to the
stochastic convolution presented in \cite{dBD1} is also given. In
section 3, the large deviation principles (LDP) is presented.
Section 4 is devoted to the support result and the two last
sections to the applications.

\section{Notations and preliminary results}\label{s2}
Throughout the paper the following notations will be used.\\
\indent For $p\in \N^*$, $\xLn^{p}(\R^d)$ is the classical
Lebesgue space of complex valued functions and $\xW^{1,p}(\R^d)$
is the associated Sobolev space of $\xLn^{p}(\R^d)$ functions with
first order derivatives, in the sense of distributions, in
$\xLn^{p}(\R^d)$. When $p=2$, $\xHn^{s}(\R^d)$ denotes the
fractional Sobolev space of tempered distributions
$v\in\mathcal{S}'(\R^d)$ such that the Fourier transform $\hat{v}$
satisfies $(1+|\xi|^2)^{s/2}\hat{v}\in \L2$. The space $\L2$ is
endowed with the inner product defined by
$(u,v)_{\L2}=\Re\int_{\R^d}u(x)\overline{v}(x)dx$. Also, when it
is clear that $\mu$ is a Borel measure on a specified Banach
space, we simply write $\xLtwo(\mu)$ and do not specify the Banach
space and Borel $\sigma-$field.\vspace{0.3cm}

\indent If $I$ is an interval of $\R$, $(E,\|\cdot\|_E)$ a Banach
space and $r$ belongs to $[1,+\infty]$, then $\xLn^{r}(I;E)$ is
the space of strongly Lebesgue measurable functions $f$ from $I$
into $E$ such that $t\rightarrow \|f(t)\|_E$ is in $\xLn^{r}(I)$.
Let $\xLn^{r}_{loc}(0,+\infty;E)$ be the respective spaces of
locally integrable functions on $(0,+\infty)$. They are endowed
with topologies of Fr\'echet space. The spaces
$\xLn^{r}(\Omega;E)$ are defined similarly.\vspace{0.3cm}

\indent We recall that a pair $(r,p)$ of positive numbers is
called an admissible pair if $p$ satisfies $2\leq
p<\frac{2d}{d-2}$ when $d>2$ ($2\leq p<+\infty$ when $d=2$ and
$2\leq p\leq+\infty$ when $d=1$) and $r$ is such that
$\frac{2}{r}=d\left(\frac{1}{2}-\frac{1}{p}\right)$. For example
$(+\infty,2)$ is an admissible pair.\vspace{0.3cm}

When $E$ is a Banach space, we will denote by $E^*$ its
topological dual space. For $x^*\in E^*$ and $x\in E$, the duality
will be denoted $<x^*,x>_{E^*,E}$.\vspace{0.3cm}

We recall that $\Phi$ is a Hilbert Schmidt operator from a Hilbert
space $H$ into a Hilbert space $\tilde{H}$ if it is a linear
continuous operator such that, given a complete orthonormal system
$(e^H_j)_{j\in\N}$ of $H$, $\sum_{j\in\N}\|\Phi
e^H_j\|^2_{\tilde{H}}<+\infty$. We will denote by
$\mathcal{L}_2(H,\tilde{H})$ the space of Hilbert Schmidt

operators from $H$ into $\tilde{H}$ endowed with the norm
$$\|\Phi\|_{\mathcal{L}_2(H,\tilde{H})}=\xtr\left(\Phi\Phi^*\right)=\sum_{j\in\N}\|\Phi
e^H_j\|_{\tilde{H}}^2,$$ where $\Phi^*$ denotes the adjoint of
$\Phi$ and $\xtr$ the trace. We denote by $\mathcal{L}_2^{s,r}$
the corresponding space for $H=\xHn^{s}(\R^d)$ and
$\tilde{H}=\xHn^{r}(\R^d)$. In the introduction $\Phi$ has been
taken in $\mathcal{L}_2^{0,1}$.\vspace{0.3cm}

\indent When $A$ and $B$ are two Banach spaces, $A\cap B$, where
the norm of an element is defined as the maximum of the norm in
$A$ and in $B$, is a Banach space. The following Banach spaces
defined for the admissible pair $(r(p),p)$ and positive $T$ by
$$X^{(T,p)}=\xCn\left([0,T];\xHone(\R^d)\right)\cap
\xL^{r(p)}\left(0,T;\xW^{1,p}(\R^d)\right)$$ will be of particular
interest.\vspace{0.3cm}

The two following Hilbert spaces of spatially localized functions
are also introduced,

$$\Sigma=\{f\in\H1:x\mapsto|x|f(x)\in\L2\}$$
endowed with the norm $\|\cdot\|_{\Sigma}$ defined by
$$\|f\|_{\Sigma}^2=\|f\|_{\H1}^2+\|x\mapsto|x|f(x)\|_{\L2}^2,$$
and $\Sigma^{\frac{1}{2}}$ where $|x|$ is replaced by
$\sqrt{|x|}$. The variance is defined as the quantity
$$V(f)=\int_{\R^d}|x|^2|f(x)|^2dx,\ f\in\Sigma.$$

Also we denote by $eval_x(f)$ the evaluation of a function $f$ at
the point $x$ where $f$ is a function taking value in any
topological space.\vspace{0.3cm}

The probability space will be denoted by
$\left(\Omega,\F,\mathbb{P}\right)$. Also, $x\wedge y$ stands for
the minimum of the two real numbers $x$ and $y$ and $x\vee y$ for
the maximum. We recall that a rate function $I$ is a lower
semicontinuous function and that a good rate function $I$ is a
rate function such that for every $c>0$, $\left\{x:I(x)\leq
c\right\}$ is a compact set. Finally, we will denote by $supp\
\mu$ the support of a probability measure $\mu$ on a topological
vector space. It is the complementary of the largest open set of
null measure.

\subsection{Properties of the group}\label{s21}
When the group acts on the Schwartz space $\mathcal{S}(\R^d)$, the
Fourier transform gives the following analytic expression
$$\forall u_0\in\mathcal{S}(\R^d),\ \forall t\neq0,\
S(t)u_0=\frac{1}{(4i\pi
t)^{\frac{d}{2}}}\int_{\R^d}e^{-i\frac{|x-y|^2}{4t}}u_0(y)dy.$$
The Fourier transform also gives that the adjoint of $S(t)$ in
$\L2$ and in every Sobolev space on $\L2$ is $S(-t)$, the same
bounded operator with time reversal.\vspace{0.3cm}

\indent The Strichartz
estimates, see \cite{kato}, are the following\\
\noindent $i/\ \ \forall\ u_0\in \L2,\ \forall\ (r,p)$ admissible
pair,
$$t\mapsto S(t)u_0 \in\xCn(\R;\L2)\cap
\xLn^{r}(\R;\xLn^{p}(\R^d)),$$ \indent$\ $ and there exists a

positive constant $c$ such that,
$$\|S(\cdot)u_0\|_{\xLn^{r}(\R;\xLn^{p}(\R^d))}\leq
c\|u_0\|_{\L2}.$$ $ii/\ $For $T>0$, for all $(r(p),p)$ and

$(r(q),q)$ two admissible pairs, if $s$ and $\rho$ are the\\
\indent$\ $conjugate exponents of $r(q)$ and $q$, i.e.
$\frac{1}{s}+\frac{1}{r(q)}=1$ and $\frac{1}{q}+\frac{1}{\rho}=1$,
$$\forall f\in
\xLn^{s}(0,T;\xLn^{\rho}(\R^d)),\ \Lambda f\in\
\xCn([0,T];\L2)\cap \xL^{r(p)}(0,T;\xLn^{p}(\R^d))$$ \indent$\
$where $\Lambda$ is defined by $\Lambda
f=\int_0^{\cdot}S(\cdot-s)f(s)ds$. Moreover, $\Lambda$ is a
continuous linear\\
\indent$\ $operator from $\xLn^{s}(0,T;\xLn^{\rho}(\R^d))$ into
$\xL^{r(p)}(0,T;\xLn^{p}(\R^d))$ with a norm that does\\
\indent$\ $not depend on $T$.
\begin{rmrk}
The first estimate gives the integrability property of the group,
the second gives the integrability of the convolution that allows
to treat the nonlinearity.
\end{rmrk}

\subsection{Topology and trajectory spaces}\label{s22}
Let us introduce a topological space that allows to treat the
subcritical case or the defocusing case. When $d>2$, we set
$$\mathcal{X}_{\infty}=\bigcap_{T\in\R^*_+,\ 2\leq p<\frac{2d}{d-2}}
X^{(T,p)},$$ it is endowed with the projective limit topology, see
\cite{BSU} and \cite{DZ}. When $d=2$ and $d=1$ we write

$p\in[2,+\infty)$.\\
The set of indices $(J,\prec)$, where $(T,p)\prec(S,q)$ if $T\leq
S$ and $p\leq q$, is a partially ordered right-filtering

set.\\
If $(T,p)\prec(S,q)$ and $u\in X^{(S,q)}$, H\"older's inequality

gives that for $\alpha$ such that
$\frac{1}{p}=\frac{\alpha}{q}+\frac{1-\alpha}{2},$
$$\exists\ c(p,q)>0:\ \|u(t)\|_{\xLn^{p}(\R^d)}\leq c(p,q)
\|u(t)\|_{\L2}^{1-\alpha}\|u(t)\|_{\xLn^{q}(\R^d)}^{\alpha}.$$
Consequently,
$$\|u(t)\|_{\xW^{1,p}(\R^d)}\leq (d+1)c(p,q)
\|u(t)\|_{\H1}^{1-\alpha}\|u(t)\|_{\xW^{1,q}(\R^d)}^{\alpha}.$$ By
time integration, along with H\"older's inequality and the fact
that $\alpha r(q)=r(p)$,
 $u$ is a function of $X^{(T,p)}$ and
\begin{equation}\label{e3} \|u\|_{X^{(T,p)}}\leq (d+1)c(p,q)
\|u\|_{X^{(S,q)}}.\end{equation} If we denote by
$p_{(S,q)}^{(T,p)}$ the dense and continuous embeddings from
$X^{(S,q)}$ into $X^{(T,p)}$, they satisfy the consistency
conditions
$$\forall\ (T,p)\prec(S,q)\prec(R,r),\
p_{(R,r)}^{(T,p)}=p_{(R,r)}^{(S,q)}\circ p_{(S,q)}^{(T,p)}.$$
Consequently, the projective limit topology is well defined by the
following neighborhood basis, given for $\varphi_1$ in

$\mathcal{X}_{\infty}$ by
$$U(\varphi_1;(T,p);\epsilon)=\left\{\varphi\in \bigcap_{(T',p')\in J}
X^{(T',p')}: \|\varphi-\varphi_1\|_{X^{(T,p)}}<\epsilon\right\}.$$
It is the weakest topology on the intersection such that for every
$(T,p)\in J$, the injection
$p_{(T,p)}:\mathcal{X}_{\infty}\rightarrow
X^{(T,p)}$ is continuous. It is a standard fact, see \cite{BSU}, that $\mathcal{X}_{\infty}$ is a Hausdorff topological space.\\
Following from (\ref{e3}), a countable neighborhood basis of
$\varphi_1$ is given by
$$\left(U\left(\varphi_1;(n,p(l));\frac{1}{k}\right)\right)_{(n,k,l)\in
(N^*)^3},$$ where $p(l)=2+\frac{4}{d-2}-\frac{1}{l}$ and
$l>\frac{d-2}{4}$ if $d>2$. If $d=2$ and $d=1$, we take $p(l)=l$.\\
It is convenient, for measurability issues, to introduce the
countable metric version $$\mathcal{D}=\bigcap_{n\in\N^*:\
n>\frac{d-2}{4}}X^{(n,p(n))};\ \ \forall(x,y)\in\mathcal{D}^2,\
d(x,y)=\sum_{n>\frac{d-2}{4}}\frac{1}{2^n}\left(\|x-y\|_{X^{(n,p(n))}}\vee
1\right).$$ It is classical that $\mathcal{D}$ is a complete
separable metric space, i.e. a Polish space. The previous
discussion gives that $\mathcal{D}$ is homeomorphic to
$\mathcal{X}_{\infty}$. Remark that it can be checked that
$\mathcal{D}$ is a locally convex Fr\'echet space.\vspace{0.3cm}

The following spaces are introduced for the case where blow-up may
occur. Adding a point $\Delta$ to the space $\H1$ and adapting
slightly the proof of Alexandroff's compactification, it can be
seen that the open sets of $\H1$ and the complementary in
$\H1\cup\{\Delta\}$ of the closed bounded sets of $\H1$ define the
open sets of a topology on $\H1\cup\{\Delta\}$. This topology
induces on $\H1$ the topology of $\H1$. Also, with such a topology
$\H1\cup\{\Delta\}$ is a Hausdorff topological space. Remark that
in \cite{AZ1}, where diffusions are studied, the compactification
of $\R^d$ is considered. Nonetheless, compactness is not an
important feature and the above construction is enough for the
following.\\ The space $\xCn([0,+\infty);\H1\cup\{\Delta\})$ is
the space of continuous functions with value in
$\H1\cup\{\Delta\}$. Also, if $f$ belongs to
$\xCn([0,+\infty);\H1\cup\{\Delta\})$ we denote the blow-up time
by
$$\mathcal{T}(f)=\inf\{t\in[0,+\infty):f(t)=\Delta\}.$$
As in \cite{AZ1}, a space of exploding paths, where $\Delta$ acts
as a cemetery, is introduced. We set
$$\mathcal{E}(\H1)=\left\{f\in\xCn([0,+\infty);\H1\cup\{\Delta\}):f(t_0)=\Delta\Rightarrow\forall
t\geq t_0,\ f(t)=\Delta\right\}.$$ It is endowed with the topology
defined by the following neighborhood basis given for $\varphi_1$
in $\mathcal{E}(\H1)$ by
$$V_{T,\ \epsilon}(\varphi_1)=\left\{\varphi\in\mathcal{E}(\H1):\mathcal{T}(\varphi)\geq \mathcal{T}(\varphi_1),\ \|\varphi_1-\varphi\|_{\xLinfty\left([0,T];\H1\right)}\leq\epsilon\right\},$$
where $T<\mathcal{T}(\varphi_1)$ and $\epsilon>0$.\\
As a consequence of the topology of $\mathcal{E}(\H1)$, the

function $\mathcal{T}$ from $\mathcal{E}(\H1)$ into $[0,+\infty]$
is sequentially lower semicontinuous, this is to say that if a
sequence of functions $(f_n)_{n\in\N}$ converges to $f$ then
$\underline{\lim}_{n\rightarrow+\infty}\mathcal{T}(f_n)\geq
\mathcal{T}(f)$. Following from (\ref{e3}), the topology of
$\mathcal{E}(\H1)$ is also defined by the countable neighborhood
basis given for $\varphi_1\in\mathcal{E}(\H1)$ by
$\left(V_{\mathcal{T}(\varphi_1)-\frac{1}{n},\frac{1}{k}}(\varphi_1)\right)_{(n,k)\in(\N^*)^2}$.
Therefore $\mathcal{T}$ is a lower semicontinuous mapping.\\
Remark that, as topological spaces, the two following spaces
satisfy the identity
$$\left\{f\in\mathcal{E}(\H1):\mathcal{T}(f)=+\infty\right\}=\xCn([0,+\infty);\H1).$$
Finally, the analogous of the intersection in the subcritical case
endowed with projective limit topology $\mathcal{E}_{\infty}$ is
defined, when $d>2$, by
$$\left\{f\in\mathcal{E}(\H1):\forall\ p\in
\left[2,\frac{2d}{d-2}\right),\ \forall\ T\in [0,\mathcal{T}(f)),\
f\in \xL^{r(p)}\left(0,T;\xW^{1,p}(\R^d)\right)\right\}.$$ When
$d=2$ and $d=1$ we write $p\in[2,+\infty)$. It is endowed with the
topology defined for $\varphi_1$ in $\mathcal{E}_{\infty}$ by the
following neighborhood basis
$$W_{T,p,\ \epsilon}(\varphi_1)=\left\{\varphi\in\mathcal{E}_{\infty}:\mathcal{T}(\varphi)\geq \mathcal{T}(\varphi_1),\ \|\varphi_1-\varphi\|_{X^{(T,p)}}\leq\epsilon\right\}.$$
where $T<\mathcal{T}(\varphi_1)$, $p$ is as above and
$\epsilon>0$. From the same arguments as for the space
$\mathcal{X}_{\infty}$, $\mathcal{E}_{\infty}$ is a Hausdorff
topological space. Also, as previously, (\ref{e3}) gives that the
topology can be defined for $\varphi_1$ in $\mathcal{E}_{\infty}$
by the countable neighborhood basis
$\left(W_{\mathcal{T}(\varphi_1)-\frac{1}{n},p(n),\frac{1}{k}}(\varphi_1)\right)_{(n,k)\in(\N^*)^2:\
n>\frac{d-2}{4}}$.\vspace{0.3cm}

If we denote again by $\mathcal{T}:\mathcal{E}_{\infty}\rightarrow
[0,+\infty]$ the blow-up time, since $\mathcal{E}_{\infty}$ is
continuously embedded into $\mathcal{E}(\H1)$, $\mathcal{T}$ is
lower semicontinuous. Thus, since the sets $\left\{[0,t],\ t\in
[0,+\infty]\right\}$ is a $\pi-$system that generates the Borel

$\s-$algebra of $[0,+\infty]$, $\mathcal{T}$ is measurable. Remark
also that, as topological spaces, the following spaces are
identical
$$\left\{f\in\mathcal{E}_{\infty}:\mathcal{T}(f)=+\infty\right\}=\mathcal{X}_{\infty}.$$

\subsection{Statistical properties of the noise} The $Q$-Wiener
process $W$ is such that its trajectories are in
$\xCn([0,+\infty);\H1)$. We assume that $Q=\Phi\Phi^*$ where
$\Phi$ is a Hilbert Schmidt operator from $\L2$ into $\H1$. The
Wiener process can therefore be written as $W=\Phi W_c$ where
$W_c$ is a cylindrical Wiener process. Recall that for any
orthonormal basis $(e_j)_{j\in\N}$ of $\L2$, there exists a
sequence of real independent Brownian motions $(\beta_j)_{j\in\N}$
such that $W_c=\sum_{j\in\N}\beta_je_j$.\vspace{0.3cm}

The sum $W_c=\sum_{j\in\N}\beta_j e_j$ is well defined in every
Hilbert space $H$ such that $\L2$ is embedded into $H$ with a
Hilbert Schmidt embedding. The denomination cylindrical is
justified by the fact that the decomposition of $W_c(1)$ on
cylinder sets $(e_1,...,e_N)$ are the finite dimensional centered
gaussian variables $(\beta_1(1),...,\beta_N(1))$ with a covariance
equal to the identity. The natural extension of the corresponding
sequence of centered gaussian measures in finite dimensions, with

a covariance equal to identity, is a gaussian cylindrical measure.
We denote it by $\nu$. The law of $W(1)$ is then the
$\sigma-$additive direct image measure
$\mu=\Phi_*\nu$.\vspace{0.3cm}

The Karhunen-Loeve decomposition of the Brownian motions on
$\xLtwo(0,1)$ is the decomposition on the eigenvectors
$\left(f_i(t)=\sqrt{2}\sin\left(\left(i+\frac{1}{2}\right)\pi t
\right)\right)_{i\in\N}$ of the injective correlation operator
$\varphi\mapsto\int_0^1(s\wedge \cdot)\varphi(s)ds$ which form a

complete orthonormal system. The coefficients are then independent
real-valued centered normal random variables where the variances
are equal to the eigenvalues $(\lambda_i)_{i\in\N}$. In addition,
$\left(g_i=\sqrt{\lambda_i}\frac{\partial}{\partial
t}f_i\right)_{i\in\N}$ also forms a complete orthonormal system of
$\xLtwo(0,1)$. Thus, formally, the coefficients of the series
expansion of the derivative of $W_c$ on the tensor product of the
complete orthonormal systems are a sequence of independent
real-valued standard normal random variables. It is thus a
gaussian white noise.\vspace{0.3cm}

The correlation operator $\tilde{Q}$ of our space-colored noise,
for $0\leq t\leq t+s\leq 1$, $(x,z)\in(\R^d)^2$,

$(h,k)\in\mathbb{C}^2$ and
$(z,t)_{\mathbb{C}}=\Re(z\overline{t})$, is formally given by
$$\E\left[\left(\frac{\p W}{\p t}(t+s,x+z),h\right)_{\mathbb{C}}\left(\frac{\p W}{\p
t}(t,x),k\right)_{\mathbb{C}}\right]=\left(\tilde{Q}k,h\right)_{\mathbb{C}},$$
where
$$\tilde{Q}k=\sum_{(i,j)\in\N^2}g_i(t+s)g_i(t)\Phi e_j(x+z)\left(\Phi
e_j(x),k\right)_{\mathbb{C}}.$$ Testing the distribution on the
product of smooth complex-valued functions with compact support on
$\R^+$ and $\R^d$, respectively $\psi$ and $\varphi$, we get
$$\int_0^1\int_{\R^d}\tilde{Q}k\psi(s)\varphi(z)dsdz=\psi(0)\overline{\Phi\Phi^*\overline{\varphi}(\cdot-x)}(x)k.$$
The correlation operator could then be identified with the
multiplication by the distribution $\delta_0\otimes
\overline{eval_x\Phi\Phi^*\overline{\tau_{-x}}}.$ Remark that if
$\Phi$ were the identity map, we would obtain the multiplication
by the Dirac mass $\delta_{(0,0)}$.\vspace{0.3cm}

The operator $\Phi$ belongs to $\mathcal{L}_2^{0,1}$ and thus to
$\mathcal{L}_2^{0,0}$ and is defined through the kernel ${\mathcal
K}(x,y)=\frac{1}{2}\sum_{j\in\N}\Phi \overline{e_j}(x)e_j(y)$ of
$\xLtwo(\R^d\times\R^d)$. This means that for any square
integrable function $u$, $\Phi u(x) = \int_{\R^d} \mathcal{K}(x,y)
u(y) dy$. Since $\Phi$ is in $\mathcal{L}_2^{0,1}$, the kernel
satisfies
$$\|\Phi\|_{\mathcal{L}_2^{0,1}}=\|\mathcal{K}\|_{\xLtwo(\R^d\times\R^d)}+\sum_{n=1}^d\sum_{j\in\N}\left\|\int_{\R^d}\frac{\partial}{\partial
x_n}\mathcal{K}(\cdot,y)e_j(y)dy\right\|_{\L2}^2.$$ Remark that
since it is impossible that $\mathcal{K}(x,y)=\mathcal{K}(x-y)$,
the noise is not homogeneous. The distribution
$\overline{eval_x\Phi\Phi^*\overline{\tau_{-x}}}$ is then the
$\L2$ function $$\frac{1}{2}\sum_{j\in\N}\Phi
e_j(x+z)\overline{\Phi e_j(x)}\ \mbox{or}\
\int_{\R^d}\mathcal{K}(x+z,u)\overline{\mathcal{K}}(x,u)du.$$

Some authors use the terminology correlation function as in our
introduction. In this setting we obtain
$$\mathbb{E}\left[\frac{\partial}{\partial
t}W(t+s,x+z)\overline{\frac{\partial}{\partial
t}W(t,x)}\right]=2\tilde{Q}$$ and
$$\mathbb{E}\left[\frac{\partial}{\partial
t}W_c(t_1,x_1)\overline{\frac{\partial}{\partial
t}W_c(t_2,x_2)}\right]=2\delta_{t_1-t_2}\otimes\delta_{x_1-x_2}$$
in the case of white noise.\vspace{0.3cm}

Finally, recall the more standard result that the correlation

operator of the $Q$-Wiener process is $(t\wedge
s)\Phi\Phi^*$.\vspace{0.3cm}

In the following we assume that the probability space is endowed

with the filtration $\mathcal{F}_t=\mathcal{N}\cup\s\{W_s,0\leq
s\leq t\}$ where $\mathcal{N}$ denotes the $\mathbb{P}-$null sets.

\subsection{The random perturbation}
Considering the fixed point problem, the stochastic convolution
$Z(t)=\int_0^tS(t-s)dW(s)$ is a stochastic perturbation term. Let
us define the operator $\mathcal{L}$ on $\xLtwo(0,T;\L2)$ by
$$\mathcal{L}h(t)= \int_0^tI\circ S(t-s)\Phi
h(s)ds,\ h\in \xLtwo(0,T;\L2),$$ where $I$ is the injection of
$\H1$ into $\L2$.

\begin{prpstn}\label{p1}
The stochastic convolution defines a measurable mapping from
$\left(\Omega,\mathcal{F}\right)$ into
$\left(\mathcal{X}_{\infty},\mathcal{B}^{X}\right)$, where
$\mathcal{B}^{X}$ stands for the Borel $\sigma-$field. Its law is
denoted by $\mu^{Z}$.\\
The direct images $\mu^{Z;(T,p)}=p_{(T,p)*}\mu^{Z}$ on the real
Banach spaces $X^{(T,p)}$ are centered gaussian measures of
reproducing kernel Hilbert space (RKHS) $H_{\mu^{Z;(T,p)}}=\xim
\mathcal{L}$ with the norm of the image structure.
\end{prpstn}
\begin{proof}Setting $F(t)=\int_0^tS(-u)dW(u)$, for $t\in\R^+$,
$Z(t)=S(t)F(t)$ follows. Indeed, if $(f_j)_{j\in\N}$ is a complete
orthonormal system of $\H1$, a straightforward calculation gives
that $(Z(t),f_j)_{\H1}=(S(t)F(t),f_j)_{\H1}$ for every $j$ in

$\N$. The continuity of the paths follows from the construction of
the stochastic integral of measurable and adapted operator
integrands satisfying
$\mathbb{E}\left[\int_0^t\|S(-u)\Phi\|_{\mathcal{L}_2^{0,0}}^2\right]<+\infty$
with respect to the Wiener process and from the strong continuity
of the group. Consequently, for every positive $T$, the paths are

in $\xCn([0,T];\H1)$.\vspace{0.3cm}

\textbf{Step 1:} The mapping $Z$ is measurable from
$\left(\Omega,\mathcal{F}\right)$ into
$\left(X^{(T,p)},\mathcal{B}^{(T,p)}\right)$,
where $\mathcal{B}^{(T,p)}$ denotes the associated Borel $\sigma-$field.\\
Since $X^{(T,p)}$ is a Polish space, every open set is a countable
union of open balls and consequently $\mathcal{B}^{(T,p)}$ is
generated by open balls.\\ Remark that the event
$\left\{\omega\in\Omega: \|Z(\omega)-x\|_{X^{(T,p)}}\leq
r\right\}$ is the intersection of
$$\bigcap_{s\in\mathbb{Q}\cap[0,T]}\left\{\omega\in\Omega:
\|Z(s)(\omega)-x\|_{\H1}\leq r\right\}$$ and of
$$\left\{\omega\in\Omega:
\|Z(\omega)-x\|_{\xL^{r(p)}(0,T;\xW^{1,p}(\R^d))}\leq r\right\}.$$
Also, remark that, since $(Z(t))_{t\in\R^+}$ is a collection of
$\H1$ random variables, the first part is a countable intersection
of elements of $\mathcal{F}$. Consequently, it suffices to show
that : $\omega\mapsto (t\mapsto Z(t))$ defines a
$\xL^{r(p)}(0,T;\xW^{1,p}(\R^d))$
random variable.\\
Consider $(\Phi_n)_{n\in\N}$ a sequence of operators of
$\mathcal{L}_2^{0,2}$ converging to $\Phi$ for the topology of
$\mathcal{L}_2^{0,1}$ and $Z_n$ the associated stochastic
convolutions. The Sobolev injections along with H\"older's
inequality give that when $d>2$ and $2\leq p\leq \frac{2d}{d-2}$,
$\H1$ is continuously embedded in $\xLn^{p}(\R^d)$. It also gives
that, when $d=2$, $\H1$ is continuously embedded in every
$\xLn^{p}(\R^d)$ for every $p\in[2,+\infty)$ and for every
$p\in[2,+\infty]$ when $d=1$. Consequently, for every $n$ in $\N$,
$Z_n$ defines a $\xCn([0,T];\xHtwo(\R^d))$ random variable and
therefore a $\xL^{r(p)}\left(0,T;\xW^{1,p}(\R^d)\right)$ random variable for the corresponding values of $p$.\\
Revisiting the proof of Proposition 3.1 in reference \cite{dBD1}
and letting $2\s+2$ be replaced by any of the previous values of
$p$ besides $p=+\infty$ when $d=1$, the necessary measurability
issues to apply the Fubini's theorem are satisfied. Also, one gets
the same estimates and that there exists a constant $C(d,p)$ such
that for every $n$ and $m$ in $\N$,
$$\mathbb{E}\left[\|Z_{n+m}(\omega)-Z_n(\omega)\|^r_{\xL^{r(p)}\left(0,T;\xW^{1,p}(\R^d)\right)}\right]\leq
C(d,p)T^{\frac{r}{2}-1}\|\Phi_{n+m}-\Phi_n\|_{\mathcal{L}_2^{0,1}}^r.$$
The sequence $(Z_n)_{n\in\N}$ is thus a Cauchy sequence of
$\xLn^{r}\left(\Omega;\xLn^{r}\left(0,T;\xW^{1,p}(\R^d)\right)\right)$,
which is a Banach space, and thus converges to $\tilde{Z}$. The
previous calculation also gives that
$$\mathbb{E}\left[\|Z_{n}(\omega)-Z(\omega)\|^r_{\xL^{r(p)}\left(0,T;\xLn^{p}(\R^d)\right)}\right]\leq
C(d,p)T^{\frac{r}{2}-1}\|\Phi_{n}-\Phi\|_{\mathcal{L}_2^{0,1}}^r.$$
Therefore $\tilde{Z}=Z$, $Z$ belongs to
$\xL^{r(p)}\left(0,T;\xW^{1,p}(\R^d)\right)$ and it defines a

measurable mapping as expected.\\
Remark that in $\mathcal{D}$, to simplify the notations, we did
not write the cases $p=+\infty$ when $d=1$ or $p=\frac{2d}{d-2}$
when $d>2$. In fact, we are interested in results on the laws of
the solutions of stochastic NLS and not really on the stochastic
convolution. Also, the result of continuity in the next section
shows that we necessarily loose on~$p$ in order to interpolate
with $2<p<p'$ and have a nonzero exponent on the $\L2-$norm.
Therefore, even if it seems at first glance that we loose on the
Sobolev's injections, it is not a restriction.\vspace{0.3cm}

\textbf{Step 2:} The mapping $Z$ is measurable with values in $\mathcal{D}$ with the Borel $\sigma-$field~$\mathcal{B}^{\mathcal{D}}$.\\
From step 1, given $x\in\mathcal{D}$, for every $n$ in
$\mathbb{N}^*$ such that $n>\frac{d-2}{4}$ the mapping

$\omega\mapsto\|Z(\omega)-x\|_{X^{(n,p(n))}}$ from
$\left(\Omega,\mathcal{F}\right)$ into
$\left(\R^+,\mathcal{B}(\R^+)\right)$, where $\mathcal{B}(\R^+)$
stands for the Borel $\s-$field of $\R^+$, is measurable. Thus
$$\omega\mapsto
d(Z(\omega),x)=\lim_{N\rightarrow+\infty}\sum_{n=1}^N\frac{1}{2^n}\left(\|Z(\omega)-x\|_{X^{(n,p(n))}}\vee
1\right)$$ is measurable. Consequently, for every $r$ in $\R^+$,
$\left\{\omega\in \Omega\
:\ d(Z(\omega),x)<r\right\}$ belongs to $\mathcal{F}$.\\
Remark that the law $\mu^{Z;\mathcal{D}}$ of $Z$ on the metric
space $\mathcal{D}$, which is a positive Borel measure, is
therefore
also regular and consequently it is a Radon measure.\\
The direct image of the Borel probability measure
$\mu^{Z;\mathcal{D}}$ by the isomorphism defines the measure
$\mu^Z$ on
$\left(\mathcal{X}_{\infty},\mathcal{B}^X\right)$.\vspace{0.3cm}

\textbf{Step 3} (Statements on the measures $\mu^{Z;(T,p)}$){\bf:}
For $(T,p)$ in the set of indices $J$, let $i_{(T,p)}$ denote the

continuous injections from $X^{(T,p)}$ into $\xLtwo(0,T;\L2)$ and
$\mu^{Z;L}=(i_{(T,p)})_*\mu^{Z;(T,p)}$. The $\sigma-$field on
$\xLtwo(0,T;\L2)$ is the Borel $\sigma-$field. Let $h\in
\xLtwo(0,T;L^2(\mathbb{R}^d))$, then
$$\left(h,i_{(T,p)}(Z)\right)_{\xLtwo(0,T;\L2)}=\int_0^T\sum_{i,j=1}^{+\infty}\int_0^t(e_j,S(t-s)\Phi
e_i)_{\L2}d\beta_i(s)(h(t),e_j)_{\L2}$$ and from classical
computation it is the almost sure limit of a sum of independent
centered gaussian random variables, thus $\mu^{Z;L}$ is a centered
gaussian measure.\\
Every linear continuous functional on $\xLtwo(0,T;\L2)$ defines by
restriction a linear continuous functional on $X^{(T,p)}$. Thus,
$\xLtwo(0,T;\L2)^*$ could be thought of as a subset of
$\left(X^{(T,p)}\right)^*$. Since $i_{(T,p)}$ is a continuous

injection, $\xLtwo(0,T;\L2)^*$ is dense in
$\left(X^{(T,p)}\right)^*$ for the weak$^*$ topology
$\s\left(\left(X^{(T,p)}\right)^*,X^{(T,p)}\right)$. It means
that, given $x^*\in\left(X^{(T,p)}\right)^*$, there exists a
sequence $\left(h_n\right)_{n\in\N}$ of elements of
$\xLtwo(0,T;\L2)$ such that for every $x\in X^{(T,p)}$,
$$\lim_{n\rightarrow+\infty}\left(h_n,i_{(T,p)}(x)\right)_{\xLtwo(0,T;\L2)}=<x^*,x>_{\left(X^{(T,p)}\right)^*,X^{(T,p)}}.$$
In other words, the random variable
$<x^*,\cdot>_{\left(X^{(T,p)}\right)^*,X^{(T,p)}}$ is a pointwise
limit of $\left(h_n,i_{(T,p)}(\cdot)\right)_{\xLtwo(0,T;\L2)}$
which are,
from the above, centered gaussian random variables. As a consequence, $\mu^{Z;(T,p)}$ is a centered gaussian measure.\\
Recall that the RKHS $H_{\mu^{Z;L}}$ of $\mu^{Z;L}$ is $\xim R^L$
where $R^L$ is the mapping from
$H_{\mu^{Z;L}}^*=\overline{\xLtwo(0,T;\L2)^*}^{\xLtwo(\mu^{Z;L})}$
with the inner product derived from the one in $\xLtwo(\mu^{Z;L})$
into $\xLtwo(0,T;\L2)$ defined for $\varphi$ in $H_{\mu^{Z,L}}^*$
by
$$R^L(\varphi)=\int_{\xLtwo(0,T;\L2)}x\varphi(x)\mu^{Z;L}(dx).$$
The same is true for $H_{\mu^{Z;(T,p)}}$ replacing
$\xLtwo(0,T;\L2)$ by $X^{(T,p)}$ and $\mu^{Z;L}$ by
$\mu^{Z;(T,p)}$.\\
Since $\mu^{Z;L}$ is the image of $\mu^{Z;(T,p)}$, taking $x^*\in
\xLtwo(0,T;\L2)^*$, we obtain that
\begin{align*}
\|x^*\|_{\xLtwo(\mu^{Z;L})}&=\int_{\xLtwo(0,T;\L2)}<x^*,x>_{\xLtwo(0,T;\L2)^*,\xLtwo(0,T;\L2)}^2\mu^{Z;L}(dx)\\
&=\int_{X^{(T,p)}}<x^*,x>_{\xLtwo(0,T;\L2)^*,\xLtwo(0,T;\L2)}^2\mu^{Z;(T,p)}(dx)\\
&=\int_{X^{(T,p)}}<x^*,x>_{\left(X^{(T,p)}\right)^*,X^{(T,p)}}^2\mu^{Z;(T,p)}(dx)=\|x^*\|_{\xLtwo(\mu^{Z;(T,p)})}.
\end{align*}
Therefore, from Lebesgue's dominated convergence theorem, we
obtain that
$$\left(X^{(T,p)}\right)^*=\overline{\xLtwo(0,T;\L2)^*}^{\s\left(\left(X^{(T,p)}\right)^*,X^{(T,p)}\right)}\subset\overline{\xLtwo(0,T;\L2)^*}^{\xLtwo(\mu^{Z;(T,p)})}$$
where the last term is equal to $H_{\mu^{Z;L}}^*$.
It follows that $H_{\mu^{Z;(T,p)}}^*\subset H_{\mu^{Z;L}}^*$.\\
The reverse inclusion follows from the fact that
$\xLtwo(0,T;\L2)^*\subset \left(X^{(T,p)}\right)^*$.\\
The conclusion follows from the quite standard fact that the RKHS
of $\mu^{Z;L}$, which is a centered gaussian measure on a Hilbert
space, is equal to $\xim \mathcal{Q}^{\frac{1}{2}}$, with the norm
of the image structure. $\mathcal{Q}$ denotes the covariance
operator of the centered gaussian measure, it is given, see
\cite{DPZ}, for $h\in \xLtwo(0,T;\L2)$, by

$$\mathcal{Q}h(v)=\int_0^T\int_0^{u\wedge v}IS(v-s)\Phi\Phi^*S(s-u)I^*h(u)dsdu.$$
Corollary B.5 of reference \cite{DPZ} finally gives that $\xim
\mathcal{L}=\xim \mathcal{Q}^{\frac{1}{2}}$.
\end{proof}

\subsection{Continuity with respect to the
perturbation}\label{s25} Recall that the mild solution of
stochastic NLS (\ref{e5}) could be written as a function of the
perturbation.\\ Let $v(x)$ denote the solution of
$$\left\{\begin{array}{cl}
    &i\frac{\xdif}{\xdif t}v-(\Delta v+|v-ix|^{2\sigma}(v-ix))=0, \\
    &v(0)=u_0,
  \end{array}\right.$$
or equivalently a fixed point of the functional
$$\mathcal{F}_x(v)(t)=S(t)u_0-i\lambda\int_0^tS(t-s)(|(v-ix)(s)|^{2\sigma}(v-ix)(s))ds,$$
where $x$ is an element of $X^{(T,p)}$, $p$ is such that
$p\geq2\s+2$ and
$(T,p)$ is an arbitrary pair in the set of indices $J$.\\
If $u$ is such that $u=v(Z)-iZ$ where $Z$ is the stochastic
convolution, note that its regularity is given in the previous
section, then $u$ is a solution of (\ref{e3}) and of (\ref{e2}).
Consequently, if $\mathcal{G}$ denotes the mapping that satisfies
$\mathcal{G}(x)=v(x)-ix$ we obtain that $u=\mathcal{G}(Z)$.\\
The local existence follows from the fact that for $R>0$ and $r>0$
fixed, taking $\|x\|_{X^{(T,2\sigma+2)}}\leq R$ and
$\|u_0\|_{\H1}\leq r$, there exists a sufficiently small
$T_{2\s+2}^*$ such that the closed ball centered at 0 of radius
$2r$ is invariant and $\mathcal{F}_x$ is a contraction for the
topology of $\xLinfty([0,T_{2\s+2}^*];\L2)\cap
\xLn^{r}(0,T_{2\s+2}^*;\xLn^{p}(\R^d))$. The closed ball is
complete for the weaker topology. The proof uses extensively the
Strichartz' estimates, see \cite{dBD1} for a detailed proof. The
same fixed point argument can be used for $\|x\|_{X^{(T,p)}}\leq
R$ in a closed ball of radius $2r$ in $X^{(T^{*}_p,p)}$ for every
$T_p^*$ sufficiently small and $p\geq 2\sigma+2$ such that
$(T_p^*,p)\in J$. From (\ref{e3}), there exists a unique maximal
solution $v(x)$ that belongs to $\mathcal{E}_{\infty}$.\\ It could
be deduced from Proposition $3.5$ of \cite{dBD1}, that the mapping
$\mathcal{G}$ from $\mathcal{X}_{\infty}$ into
$\mathcal{E}_{\infty}$ is a continuous mapping from
$\bigcap_{T\in\R^*_+}X^{(T,2\sigma+2)}$ with the projective limit
topology into $\mathcal{E}(\H1)$. The result can be strengthen as
follows.
\begin{prpstn}\label{p2}
The mapping $\mathcal{G}$ from $\mathcal{X}_{\infty}$ into
$\mathcal{E}_{\infty}$ is continuous.
\end{prpstn}
\begin{proof} Let $\tilde{x}$ be a function of $\mathcal{X}_{\infty}$ and
$T<\mathcal{T}(\tilde{Z})$. Revisiting the proof of Proposition
$3.5$ of \cite{dBD1} and taking $\epsilon>0$, $p'\geq 2\sigma+2$,
$2<p<p'$, $R=1+\|\tilde{x}\|_{X^{(T,p')}}$ and
$r=1+\|v(\tilde{x})\|_{\xCn([0,T];\H1)}$ there exists $\eta$
satisfying
$$0<\eta<\frac{\epsilon}{2(d+1)C(p,p')}\wedge 1$$ such for
$x$ in $\mathcal{X}_{\infty}$
$$\|x-\tilde{x}\|_{X^{(T,p')}}\leq\eta\Rightarrow\|v(x)-v(\tilde{x})\|_{\xCn([0,T];\H1)}\leq\left(\frac{\epsilon}{2(d+1)C(p,p')(4r)^{\alpha}}\right)^{\frac{1}{1-\alpha}}\wedge
1.$$ The constant $\alpha$ is the one that appears in the
application of H\"older's inequality before (\ref{e3}).
Consequently, since $v(x)$ and $v(\tilde{x})$ are functions of the
closed ball centered at 0 and of radius $2r$ in $X^{(T,p)}$, the
triangular inequality gives that
$$\|v(x)-v(\tilde{x})\|_{X^{(T,p')}}\leq 4r.$$The application of
both H\"older's inequality and the triangular inequality allow to
conclude that
$$\forall
x\in\mathcal{X}_{\infty}:\|x-\tilde{x}\|_{X^{(T,p')}}\leq\eta,\
\|\mathcal{G}(x)-\mathcal{G}(\tilde{x})\|_{X^{(T,p)}}\leq\epsilon$$
which, from the definition of the neighborhood basis of
$\mathcal{E}_{\infty}$, gives the continuity. \end{proof}

The following corollary is a consequence of the last statement of

section \ref{s22}.
\begin{crllr}
In the focusing subcritical case or in the defocusing case,
$\mathcal{G}$ is a continuous mapping from $\mathcal{X}_{\infty}$
into $\mathcal{X}_{\infty}$
\end{crllr}
The continuity allows to define the law of the solutions of the
stochastic NLS equations on $\mathcal{E}_{\infty}$ and in the
cases of global existence in $\mathcal{X}_{\infty}$ as the direct
image $\mu^u=\mathcal{G}_*\mu^Z$, the same
notation will be used in both cases.\\
Let consider the solutions of
\begin{equation}
\label{e4} i \xdif u_{\epsilon} - (\Delta u_{\epsilon} + \lambda
|u_{\epsilon}|^{2\s}u_{\epsilon})\xdif t = \sqrt{\epsilon}\xdif W,
\end{equation}
where $\epsilon\geq 0$. The laws of the solutions $u_{\epsilon}$
in the corresponding trajectory spaces are denoted by
$\mu^{u_{\epsilon}}$, or equivalently
$\mathcal{G}_*\mu^{Z_{\epsilon}}$ where $\mu^{Z_{\epsilon}}$ is
the direct image of $\mu^Z$ under the transformation
$x\mapsto\sqrt{\epsilon}x$ on $\mathcal{X}_{\infty}$. The
continuity also gives that the family converges weakly to the
Dirac mass on the deterministic solution $u_d$ as $\epsilon$
converges to zero. Next section is devoted to the study of the
convergence towards $0$ of rare events or tail events of the law
of the solution $u_{\epsilon}$, namely large deviations. It allows
to describe more precisely the convergence towards the
deterministic measure.

\section{Sample path large deviations}\label{s3}
\begin{thrm}
The family of probability measures

$\left(\mu^{u_{\epsilon}}\right)_{\epsilon\geq0}$ on
$\mathcal{E}_{\infty}$ satisfies a LDP of speed $\epsilon$ and
good rate function
$$I(u)=\frac{1}{2}\inf_{h\in\xLtwo(0,+\infty;\L2):\mathbf{S}(h)=u}\left\{\|h\|_{\xLtwo(0,+\infty;\L2)}^2\right\},$$
where $\inf\emptyset=+\infty$ and $\mathbf{S}(h)$, called the
skeleton, is the unique mild solution of the following control
problem:
$$\left\{\begin{array}{cl}
&i\frac{\xdif}{\xdif t}u=\Delta u +\lambda |u|^{2\sigma}u + \Phi h,\\
&u(0)=u_0\ \in \xHone(\R^d).\\
\end{array}\right.$$ 
This is to say that for every Borel set $A$ of
$\mathcal{E}_{\infty}$,
$$-\inf_{u\in\AA}I(u)\leq\underline{\lim}_{\epsilon\rightarrow 0}
\epsilon\log\mu^{u_{\epsilon}}(A)\leq
\overline{\lim}_{\epsilon\rightarrow
0}\epsilon\log\mu^{u_{\epsilon}}(A)\leq-\inf_{u\in\overline{A}}I(u).$$ 
The same
result holds in $\mathcal{X}_{\infty}$ for the family of laws of
the solutions in the cases of global existence.
\end{thrm}
\begin{proof} The general LDP for centered Gaussian measures on real
Banach spaces, see \cite{DS1}, gives that for a given pair $(T,p)$
in the set of indices $J$, the family
$\left(p_{(T,p)*}\mu^{Z_{\epsilon}}\right)_{\epsilon\geq0}$
satisfies a LDP on $X^{(T,p)}$ of speed $\epsilon$ and good rate
function defined for $x\in X^{(T,p)}$ by,
$$I^{Z;(T,p)}(x)=\left\{\begin{array}{cl}
&\frac{1}{2}\|x\|_{H_{\mu^{Z;(T,p)}}}^2,\ \ \ \mbox{if}\ x\in H_{\mu^{Z;(T,p)}},\\
&+\infty,\hspace{1.8cm} \mbox{otherwise},\\
\end{array}\right.$$
which, using Proposition \ref{p1}, is equal to
$$I^{Z;(T,p)}(x)=\frac{1}{2}\inf_{h\in\xLtwo(0,T;\L2):\mathcal{L}(h)=x}\left\{\|h\|_{\xLtwo(0,T;\L2)}^2\right\}.$$ Dawson-G\"artner's theorem, see \cite{DZ}, allows to deduce that the
family $\left(\mu^{Z_{\epsilon}}\right)_{\epsilon\geq0}$ satisfies

the LDP with the good rate function defined for
$x\in\mathcal{X}_{\infty}$ by
\begin{align*}I^Z(x)&=\sup_{(T,p)\in
J}\left\{I^{Z;(T,p)}(x)\right\}\\&=\frac{1}{2}\inf_{h\in\xLtwo(0,T;\L2):\mathcal{L}(h)=x}\left\{\|h\|_{\xLtwo(0,+\infty;\L2)}^2\right\}.\end{align*}
It has been shown in sections \ref{s22} and \ref{s25} that
$\mathcal{G}$ is a continuous function from a Hausdorff
topological space into another Hausdorff topological space.
Consequently, both results follow from  Varadhan's contraction
principle along with the fact that if
$\mathcal{G}\circ\mathcal{L}(h)=x$ then $x$ is the unique mild
solution of the control problem (i.e. $x=\mathbf{S}(h)$).
\end{proof}
\begin{rmrk}
The rate function is such that
$$I(u)=\frac{1}{2}\int_0^{\mathcal{T}(u)}\left\|\left(\Phi_{|\xker
\Phi^{\bot}}\right)^{-1}\left(i\frac{\xdif }{\xdif t}u-\Delta
u-\lambda|u|^{2\sigma}u\right)(s)\right\|_{\L2}ds,$$if
$i\frac{\xdif}{\xdif t} u-\Delta u-\lambda|u|^{2\sigma}u \in \xim
\Phi$, and $I(u)=+\infty$ otherwise.
\end{rmrk}
\begin{rmrk}
In the cases where blow-up may occur, the argument that will
follow allows to prove the weaker result that, given an $(T,p)$ in
the set of indices $J$ and
$$I^{(T,p)}(u)=\frac{1}{2}\inf_{h\in\xLtwo(0,T;\L2):S(h)=u}\left\{\|h\|_{\xLtwo(0,T;\L2)}^2\right\},$$
then for every bounded Borel set $A$ of $X^{(T,p)}$
$$-\inf_{u\in \AA}I^{(T,p)}(u)\leq\underline{\lim}_{\epsilon\rightarrow0}\epsilon\log\mathbb{P}\left(u_{\epsilon}\in A\right)
\leq\overline{\lim}_{\epsilon\rightarrow0}\epsilon\log\mathbb{P}\left(u_{\epsilon}\in
A\right)\leq-\inf_{u\in \overline{A}}I^{(T,p)}(u).$$
\end{rmrk}
\noindent Indeed, if $u_{\epsilon}$ belongs to $A$, there exists a
constant $R$ such that $\|u_{\epsilon}\|_{X^{(T,p)}}\leq R$.
Denoting by $u_{\epsilon}^R$ the solution of the following fixed
point problem
$$u_{\epsilon}^R(t)=S(t)u_0-i\lambda\int_0^tS(t-s)(|(u_{\epsilon}^R-i\sqrt{\epsilon}Z)(s)|^{2\sigma}(u_{\epsilon}^R-i\sqrt{\epsilon}Z)(s))\indic_{\|u_{\epsilon}^R\|_{X^{(s,p)}}\leq R}ds,$$
the arguments used previously allow to show that
$\sqrt{\epsilon}Z\rightarrow u_{\epsilon}^R$ is a continuous
mapping from every $X^{(T,p')}$ into $X^{(T,p)}$ for $p'>p$. The
result on the laws of $u_{\epsilon}^R$ follows from Varadhan's
contraction principle replacing $S(h)$ by $S^R(h)$ with the
truncation in front of the nonlinearity. Finally, the statement
follows from the fact that $\|u_{\epsilon}\|_{X^{(T,p)}}\leq R$
implies that $u_{\epsilon}^R=u_{\epsilon}$ and that
$$\inf_{h\in\xLtwo(0,T;\L2):S^R(h)\in\overline{A}}\{\|h\|_{\xLtwo(0,T;\L2)}^2\}=\inf_{h\in\xLtwo(0,T;\L2):S(h)\in\overline{A}}\{\|h\|_{\xLtwo(0,T;\L2)}^2\}.$$

Remark that writing $\frac{\partial}{\partial t}h$ instead of $h$
in the optimal control problem leads to a rate function consisting
in the minimisation of
$\frac{1}{2}\|h\|^2_{H_0^1(0,+\infty;\L2)}$. This space is somehow
the equivalent of the Cameron-Martin space for the Brownian
motion. Specifying only the law $\mu$ of $W(1)$ on $\H1$ and

dropping $\Phi$ in the control problem would lead to a rate
function consisting in the minimisation of
$\frac{1}{2}\|h\|^2_{H_0^1(0,+\infty;H_{\mu})}$, where $H_{\mu}$
stands for the RKHS of $\mu$.\vspace{0.3cm}

The formalism of a LDP stated in the intersection space with a
projective limit topology allows, for example, to deduce by
contraction, when there is no blow-up in finite time, a variety of
sample path LDP on every $X^{(T,p)}$. The rate function could be
interpreted as the minimal energy to implement
control.\vspace{0.3cm}

LDP for the family of laws of $u_{\epsilon}(T)$, for a fixed $T$,
could be deduced by contraction in the cases of global existence.
The rate function is then the minimal energy needed to transfer
$u_0$ to $x$ from 0 to $T$. An application of this type will be
given in section \ref{s6}.\vspace{0.3cm}

Next section gives a characterization of the support of the law of
the solution in our setting. Section \ref{s5} is devoted to some
consequences of these results on the blow-up time. Finally, in
section \ref{s6}, applications in nonlinear optics are given.

\section{Remark on the support of the law of the
solution}\label{s4}
\begin{thrm}[The support theorem]
The support of the law of the solution is characterized by
$$supp\ \mu^u=\overline{\xim \mathbf{S}}^{\mathcal{E}_{\infty}}$$
and in the cases of global existence by
$$supp\ \mu^u=\overline{\xim \mathbf{S}}^{\mathcal{X}_{\infty}}$$
\end{thrm}
\begin{proof} \textbf{Step 1: }From Proposition \ref{p2}, given $(T,p)$ in the set of indices $J$, $\mu^{Z;(T,p)}$ is a gaussian measure on a Banach
space and its RKHS is $\xim \mathcal{L}$. Consequently, see
\cite{BC} Theorem (IX,2;1), its support is $\overline{\xim
\mathcal{L}}^{X^{(T,p)}}$. Also, from the definition of the image
measure we have that
$$\mu^Z\left(p_{(T,p)}^{-1}\left(\overline{\xim
\mathcal{L}}^{X^{(T,p)}}\right)\right)=\mu^{Z;(T,p)}\left(\overline{\xim
\mathcal{L}}^{X^{(T,p)}}\right)=1.$$ As a consequence the first
inclusion follows
$$supp\
\mu^Z\subset\bigcap_{(T,p)}p_{(T,p)}^{-1}\left(\overline{\xim
\mathcal{L}}^{X^{(T,p)}}\right)=\overline{\xim
\mathcal{L}}^{\mathcal{X}_{\infty}}.$$ It then suffices to show
that $\xim \mathcal{L}\subset supp\ \mu^Z$. Suppose that $x\notin
supp\ \mu^Z$, then there exists a neighborhood $V$ of $x$ in
$\mathcal{X}_{\infty}$, satisfying
$V=\bigcap_{i=1}^nV^{(T_i,p_i)}$ where $V^{(T_i,p_i)}$ is a
neighborhood of $x$ in $X^{(T_i,p_i)}$, $n$ is a finite integer
and $(T_i,p_i)$ a finite sequence of elements of $J$, such that
$\mu^Z(V)=0$. It can be shown that $\bigcap_{i=1}^n X^{(T_i,p_i)}$
is still a separable Banach space. It is such that $\mathcal{D}$
is continuously embedded into it, and such that the Borel direct
image probability measure is a Gaussian measure of RKHS $\xim
\mathcal{L}$. The support of this measure is then the adherence of
$\xim \mathcal{L}$ for the topology defined by the maximum of the
norms on each factor. Thus, $V\cap
\xim \mathcal{L}=\emptyset$ and $x\notin \xim \mathcal{L}$.\\
\textbf{Step 2:}  We conclude using the continuity of
$\mathcal{G}$.\\
Indeed since $\mathcal{G}(\xim
\mathcal{L})\subset\overline{\mathcal{G}(\xim
\mathcal{L})}^{\mathcal{E}_{\infty}}$, $\xim
\mathcal{L}\subset\mathcal{G}^{-1}\left(\overline{\mathcal{G}(\xim
\mathcal{L})}^{\mathcal{E}_{\infty}}\right)$. Since $\mathcal{G}$
is continuous, the right side is a closed set of
$\mathcal{X}_{\infty}$ and from step 1,
$$supp\ \mu^Z\subset\mathcal{G}^{-1}\left(\overline{\xim\left(
\mathcal{G}\circ\mathcal{L}\right)}^{\mathcal{E}_{\infty}}\right),$$
and
$$\mu^Z\left(\mathcal{G}^{-1}\left(\overline{\xim
\mathbf{S}}^{\mathcal{E}_{\infty}}\right)\right)=1,$$ thus
$$supp\ \mu^u\subset\overline{\xim \mathbf{S}}^{\mathcal{E}_{\infty}}.$$
Suppose that $x\notin supp\ \mu^u$, there exists a neighborhood
$V$ of $x$ in $\mathcal{E}_{\infty}$ such that
$\mu^u(V)=\mu^Z\left(\mathcal{G}^{-1}(V)\right)=0$, consequently
$\mathcal{G}^{-1}(V)\bigcap \xim \mathcal{L}=\emptyset$ and
$x\notin \xim \mathbf{S}$. This gives reverse inclusion.\\
The same arguments hold replacing $\mathcal{E}_{\infty}$ by
$\mathcal{X}_{\infty}$.
\end{proof}

Remark that the result of step 2 is general and gives  that the
support of the direct images $\mu^E$ of the law $\mu^u$ by any
continuous mapping $f$ from either $\mathcal{E}_{\infty}$ or
$\mathcal{X}_{\infty}$ into a topological vector space $E$ is
$\overline{\xim \left(f\circ \mathbf{S}\right)}^{E}$. For example,
in the cases of global existence, given a positive $T$, the
support of the law in $\H1$ of $u(T)$ is $\overline{\xim
\mathbf{S}(T)}^{\H1}$.

\section{Applications to the blow-up time}\label{s5}
In this section the equation with a focusing nonlinearity, i.e.
$\lambda=1$, is considered. In this case, it is known that some
solutions of the deterministic equation blow up in finite time. It
has been proved in section \ref{s22} that $\mathcal{T}$ is a
measurable mapping from $\mathcal{E}_{\infty}$ to $[0,+\infty]$,
both spaces are equipped with their Borel $\s-$fields.
Incidentally, $\mathcal{T}(u)$ is a $\mathcal{F}_t-$stopping time.
Also, if $B$ is a Borel set of $[0,+\infty]$,
$$\mathbb{P}\left(\mathcal{T}(u)\in
B\right)=\mu^{u}\left(\mathcal{T}^{-1}(B)\right).$$

The support theorem allows to determine whether an open or a
closed set of the form $\mathcal{T}^{-1}(B)$ is such that
$\mu^{u}\left(\mathcal{T}^{-1}(B)\right)>0$ or
$\mu^{u}\left(\mathcal{T}^{-1}(B)\right)<1$ respectively. An
application of this fact is given in Proposition \ref{p11}. For a
Borel set $B$ such that the interior of
$\mathcal{T}^{-1}(B)\cap\overline{\xim
\mathbf{S}}^{\mathcal{E}_{\infty}}$ is nonempty,
$\mathbb{P}(\mathcal{T}(u)\in B)>0$ holds. \vspace{0.3cm}

Also, $\mathcal{T}$ is not continuous and Varadhan's contraction
principle does not allow to obtain a LDP for the law of the
blow-up time. Nonetheless, the LDP for the family
$\left(\mu^{u_{\epsilon}}\right)_{\epsilon>0}$ gives the
interesting result that
$$-\inf_{u\in
Int\left(\mathcal{T}^{-1}(B)\right)}I(u)\leq\underline{\lim}_{\epsilon\rightarrow0}\epsilon\log\mathbb{P}\left(\mathcal{T}(u_{\epsilon})\in
B\right)$$and that
$$\overline{\lim}_{\epsilon\rightarrow0}\epsilon\log\mathbb{P}\left(\mathcal{T}(u_{\epsilon})\in
B\right)\leq-\inf_{u\in \overline{\mathcal{T}^{-1}(B)}}I(u),
$$
where $Int\left(\mathcal{T}^{-1}(B)\right)$ stands for the
interior set of $\mathcal{T}^{-1}(B)$. Remark also that the
interior or the adherence of sets in $\mathcal{E}_{\infty}$ are
not really tractable. In that respect, the semicontinuity of
$\mathcal{T}$ makes the sets $(T,+\infty]$ and $[0,T]$
particularly interesting.

\subsection{Probability of blow-up after time T}\label{s51}
\begin{prpstn}\label{p11}
If $u_0\in \xHn^{3}(\R^d)$ and $\xker \Phi^*=\{0\}$ then for every
positive $T$,
$$\mathbb{P}(\mathcal{T}(u)>T)>0.$$
\end{prpstn}
\begin{proof}Since $\mathcal{T}$ is lower semicontinuous,
$\mathcal{T}^{-1}((T,+\infty])$ is an open set.\\
Consider $H=-\Delta u_0-|u_0|^{2\s}u_0$ which satisfies
$\mathcal{G}\circ\Lambda(H)=u_0$, where $\Lambda$ has been defined
in section \ref{s21}, then $\mathcal{T}(\mathbf{S}(H))=+\infty$.
Also, $\Phi$ defines an operator from
$\xLtwo_{loc}(0,+\infty;\L2)$ into $\xLtwo_{loc}(0,+\infty;\H1)$
and it can be shown, since $\xker \Phi^*=\{0\}$, that its range is
dense. Consequently, there exists a sequence $(h_n)_{n\in\N}$ of
$\xLtwo_{loc}(0,+\infty;\L2)$ functions such that
$(\Phi(h_n))_{n\in\N}$
converges to $H$ in $\xLtwo_{loc}(0,+\infty;\H1)$.\\
Using the semicontinuity of $\mathcal{T}$, the continuity of
$\mathcal{G}$, the fact that
$\mathbf{S}=\mathcal{G}\circ\Lambda\circ\Phi$, the following Lemma
and the fact that $\xLtwo_{loc}(0,+\infty;\H1)$ is continuously
embedded in $\xLone_{loc}(0,+\infty;\H1)$,
$\underline{\lim}_{n\rightarrow\infty}\mathcal{T}(\mathbf{S}(h_n))\geq+\infty$,
i.e.
$\lim_{n\rightarrow\infty}\mathcal{T}(\mathbf{S}(h_n))=+\infty$,
follows.
Therefore $\mathcal{T}(\mathbf{S}(h_n))>T$ for $n$ large enough and $\mathcal{T}^{-1}((T,+\infty])\cap(\xim \mathbf{S})$ is nonempty.\\
The conclusion follows then from the support theorem.
\end{proof}
As a corollary, taking the complementary of
$\mathcal{T}^{-1}((T,+\infty])$, $\mathbb{P}(\mathcal{T}(u)\leq
T)<1$ follows. This is related to the results of \cite{dBD2} where
it is proved that for every positive $T$,
$\mathbb{P}(\mathcal{T}(u)<T)>0$ and to the graphs in section $4$
of \cite{DMD}.
\begin{lmm}
The operator $\Lambda$ from $\xLone_{loc}(0,+\infty;\H1)$ into
$\mathcal{X}_{\infty}$ of $ii/$ of the Strichartz estimates is
continuous.
\end{lmm}
\begin{proof}
The result follows from $ii/$ of the Strichartz estimates, the
fact that the partial derivatives with respect to one space
variable commutes with both the integral and the group and the
definition of the projective limit topology.
\end{proof}\vspace{0.3cm}

The following result holds when the intensity of the noise
converges to zero.

\begin{prpstn}\label{p3} If $u_0\in \xHn^{3}(\R^d),\ \xker \Phi^*=\{0\}$ and
$T\geq\mathcal{T}(u_d)$, where $u_d$ is the solution of the
deterministic NLS equation with initial datum $u_0$, there exists
a positive constant c such that
$$\underline{\lim}_{\epsilon\rightarrow0}\epsilon\log\mathbb{P}\left(\mathcal{T}(u_{\epsilon})>T\right)\geq
-c.$$
\end{prpstn}
\begin{proof}The result follows from
$$-\frac{1}{2}\inf_{h\in\xLtwo(0,+\infty;\L2):\mathcal{T}(\mathbf{S}(h))>T}\left\{\|h\|_{\xLtwo(0,+\infty;\L2)}^2\right\}\leq\underline{\lim}_{\epsilon\rightarrow0}\epsilon\log\mathbb{P}\left(\mathcal{T}(u_{\epsilon})>T\right)$$
and the fact that, from the arguments of the proof of Proposition
\ref{p11}, for every $T$ such that $T\geq\mathcal{T}(u_d)$ the set
$\{h\in \xLtwo(0,+\infty;\L2):\mathcal{T}(\mathbf{S}(h)>T\}$ is
nonempty.
\end{proof}
In the following we will denote by $L^{(T,+\infty)}$ the infimum
in the left hand side of the inequality of the above proof.
\begin{rmrk}
The assumption that $u_0\in \xHn^{3}(\R^d)$ could be dropped using
similar arguments as in Proposition $3.3$ of \cite{dBD2}.
\end{rmrk}
Remark that the LDP does not give interesting information on the
upper bound even if the bounds have been sharpened using the

rather strong projective limit topology. It is zero since $h=0$
belongs to $\overline{\mathcal{T}^{-1}((T,+\infty])}$ as for every
$T>0$,
$\overline{\mathcal{T}^{-1}((T,+\infty])}=\mathcal{E}_{\infty}$.
Indeed, if a function $f$ of $\mathcal{E}_{\infty}$ is given and
blows up at a particular time $\mathcal{T}(f)$ such that
$T>\mathcal{T}(f)$, it is possible to build a sequence
$(f_n)_{n\in \N}$ of functions of $\mathcal{E}_{\infty}$ equal to
$f$ on $\left[0,\mathcal{T}(f)-\frac{1}{n}\right]$ and such that
$\mathcal{T}(f_n)>T$. The same problem will appear in the next
section where the LDP gives a lower bound equal to $-\infty$.
Indeed, $Int\left(\mathcal{T}^{-1}([0,T])\right)$ is the
complementary of the above and thus an empty set. To overcome this
problem the approximate blow-up time is introduced. Remark also
that it is possible that $L^{(T,+\infty)}=0$.\vspace{0.3cm}

Also, the case $T<\mathcal{T}(u_d)$ has not been treated. Indeed,
the associated event is not a large deviation event and the LDP
only gives that
$$\lim_{\epsilon\rightarrow0}\epsilon\log\mathbb{P}\left(\mathcal{T}(u_{\epsilon})>
T\right)=0.$$

\subsection{Probability of blow-up before time T}
In that case we obtain
$$-\infty\leq\underline{\lim}_{\epsilon\rightarrow0}\epsilon\log\mathbb{P}\left(\mathcal{T}(u_{\epsilon})\leq T\right)\leq
\overline{\lim}_{\epsilon\rightarrow0}\epsilon\log\mathbb{P}\left(\mathcal{T}(u_{\epsilon})\leq
T\right)\leq-U^{[0,T]}$$ where
$U^{[0,T]}=\frac{1}{2}\inf_{h\in\xLtwo(0,+\infty;\L2):\mathcal{T}(\mathbf{S}(h))\leq
T}\left\{\|h\|_{\xLtwo(0,+\infty;\L2)}^2\right\}$.
\begin{prpstn}\label{p4}
If $T<\mathcal{T}(u_d)$,
$$\overline{\lim}_{\epsilon\rightarrow0}\epsilon\log\mathbb{P}\left(\mathcal{T}(u_{\epsilon})\leq
T\right)\leq-U^{[0,T]}<0.$$ Moreover, if $u_0\in\xHn^{3}(\R^d)$
and if $u_0$, $\Delta u_0$ and $|u_0|^{2\s}u_0$ belong to $\xim
\Phi$ then
$$-\infty<\underline{\lim}_{\epsilon\rightarrow0}\epsilon\log\mathbb{P}\left(\mathcal{T}(u_{\epsilon})\leq
T\right).$$

\end{prpstn}
\begin{proof} Let $(h_n)_{n\in\N}$ be a sequence of
$\xLtwo(0,+\infty;\L2)$ functions converging to zero. It follows
from Lemma 5.2 and the fact that $\xLtwo(0,+\infty;\H1)$ is
continuously embedded into $\xLone_{loc}(0,+\infty;\H1)$ that
$\mathbf{S}=\mathcal{G}\circ\Lambda\circ\Phi$ is continuous from
$\xLtwo(0,+\infty;\L2)$ into $\mathcal{E}_{\infty}$.\\ Also, from
the semicontinuity of $\mathcal{T}$,
$\underline{\lim}_{\epsilon\rightarrow0}\mathcal{T}(\mathbf{S}(h_n))\geq\mathcal{T}(u_d)$
and the first point follows.\\
The $\xLtwo(0,+\infty;\L2)$ control
\begin{align*}
H^{\mathcal{E}}(t)=&\frac{2}{T-2t}\indic_{t\leq
\frac{T}{2}}\left[-\frac{2i}{T-2t}\left(\Phi_{|\xker
\Phi^{\bot}}\right)^{-1}u_0-\left(\Phi_{|\xker
\Phi^{\bot}}\right)^{-1}(\Delta u_0)\right.\\
&\left.-\left(\frac{2}{T-2t}\right)^2\left(\Phi_{|\xker
\Phi^{\bot}}\right)^{-1}(|u_0|^{2\s}u_0)\right]
\end{align*} is
such that $\mathbf{S}(H^{\mathcal{E}})=\frac{2}{T-2t}u_0$ which
blows up before $T$. This proves the second point.\end{proof} When
$T\geq \mathcal{T}(u_d)$, the probability is not supposed to tend
to zero. Also, as $h=0$ is a solution, the upper bound is zero and
none of the bounds are interesting.

\subsection{Bounds for the approximate blow-up time}
To overcome the limitation that
$\overline{\mathcal{T}^{-1}((T,+\infty))}=\mathcal{E}_{\infty}$,
which does not allow to have two interesting bounds
simultaneously, we introduce for every positive $R$ the mappings
$\mathcal{T}_R$ defined for $f\in\mathcal{E}_{\infty}$ by
$$\mathcal{T}_R(f)=\inf\{t\in[0,+\infty):\ \|f(t)\|_{\H1}\geq R\}.$$ It
corresponds to the approximation of the blow-up time used in
\cite{DMD}. We obtain the following bounds.
\begin{prpstn}\label{p5} When $T\geq\mathcal{T}_R(u_d)$, the
following inequality holds
$$-c<-L_R^{(T,+\infty)}\leq\underline{\lim}_{\epsilon\rightarrow0}\epsilon\log\mathbb{P}\left(\mathcal{T}_R(u_{\epsilon})>T\right)$$and$$\overline{\lim}_{\epsilon\rightarrow0}\epsilon\log\mathbb{P}\left(\mathcal{T}_R(u_{\epsilon})>T\right)\leq-\sup_{\alpha>0}L_{R+\alpha}^{(T,+\infty)}.$$
Also, when $T<\mathcal{T}_R(u_d)$, we have that
$$-\inf_{\alpha>0}U_{R+\alpha}^{[0,T]}\leq\underline{\lim}_{\epsilon\rightarrow0}\epsilon\log\mathbb{P}\left(\mathcal{T}_R(u_{\epsilon})\leq
T\right)$$and$$
\overline{\lim}_{\epsilon\rightarrow0}\epsilon\log\mathbb{P}\left(\mathcal{T}_R(u_{\epsilon})\leq
T\right)\leq-U_R^{[0,T]}<0.$$In the above, $c$ is nonnegative and
the numbers $L_R^{(T,+\infty)}$ and $U_R^{[0,T]}$ are defined as
$L^{(T,+\infty)}$ and $U^{[0,T]}$ replacing $\mathcal{T}$ by
$\mathcal{T}_R$.
\end{prpstn}
\begin{proof} The result follows from the facts that $\mathcal{T}_R$,
which is not continuous, is lower semicontinuous, that for every
positive $\alpha$,
$\overline{\mathcal{T}_R^{-1}((T,+\infty))}\subset\mathcal{T}_{R+\alpha}^{-1}((T,+\infty))$,

thus $\mathcal{T}_{R+\alpha}^{-1}([0,T])\subset
Int\left(\mathcal{T}_R^{-1}([0,T])\right)$ and from the arguments
used in the proofs of the last two propositions.
\end{proof}
We also obtain the following estimates of other large deviation
events.
\begin{crllr}
If $S,T<\mathcal{T}_R(u_d)$, for every positive $c$, there exists
a positive $\epsilon_0$ such that if $\epsilon\leq \epsilon_0$,
$$L_{<,R,\epsilon,c}^{S,T}\leq\mathbb{P}\left(S<\mathcal{T}_R(u_{\epsilon})\leq T\right)\leq U_{<,R,\epsilon,c}^{S,T}$$
where
$$L_{<,R,\epsilon,c}^{S,T}=\exp\left(-\frac{\inf_{\alpha>0}U_{R+\alpha}^{[0,T]}+c}{\epsilon}\right)\left(1-\exp\left(-\frac{U_R^{[0,S]}-\inf_{\alpha>0}U_{R+\alpha}^{[0,T]}}{\epsilon}\right)\right)$$

and
$$U_{<,R,\epsilon,c}^{S,T}=\exp\left(-\frac{U_R^{[0,T]}-c}{\epsilon}\right)\left(1-\exp\left(-\frac{\inf_{\alpha>0}U_{R+\alpha}^{[0,S]}-U_R^{[0,T]}}{\epsilon}\right)\right).$$
Also, if $S,T>\mathcal{T}_R(u_d)$,  for every positive $c$, there
exists a positive $\epsilon_0$ such that if $\epsilon\leq
\epsilon_0$,
$$L_{>,R,\epsilon,c}^{S,T}\leq\mathbb{P}\left(S<\mathcal{T}_R(u_{\epsilon})\leq
T\right)\leq U_{>,R,\epsilon,c}^{S,T}$$ where
$$L_{>,R,\epsilon,c}^{S,T}=\exp\left(-\frac{L_R^{(S,+\infty)}+c}{\epsilon}\right)\left(1-\exp\left(-\frac{\sup_{\alpha>0}L_{R+\alpha}^{(T,+\infty)}-L_R^{(T,+\infty)}}{\epsilon}\right)\right)$$
and
$$U_{>,R,\epsilon,c}^{S,T}=\exp\left(-\frac{\sup_{\alpha>0}L_{R+\alpha}^{(S,+\infty)}-c}{\epsilon}\right)\left(1-\exp\left(-\frac{L_R^{(T,+\infty)}-\sup_{\alpha>0}L_{R+\alpha}^{(S,+\infty)}}{\epsilon}\right)\right).$$
\end{crllr}
\begin{proof} When $S,T<\mathcal{T}_R(u_d)$, the result
follows from the inequalities and from the fact that
\begin{align*}
\mathbb{P}\left(S<\mathcal{T}_R(u_{\epsilon})\leq
T\right)&=\mathbb{P}\left(\{\mathcal{T}_R(u_{\epsilon})\leq
T\}\setminus\{\mathcal{T}_R(u_{\epsilon})\leq
S\}\right)\\&=\mathbb{P}\left(\mathcal{T}_R(u_{\epsilon})\leq
T\right)
\left(1-\frac{\mathbb{P}\left(\mathcal{T}_R(u_{\epsilon})\leq
S\right)}{\mathbb{P}\left(\mathcal{T}_R(u_{\epsilon})\leq
T\right)}\right).\end{align*} When $S,T>\mathcal{T}_R(u_d)$, we
use
\begin{align*}\mathbb{P}\left(S<\mathcal{T}_R(u_{\epsilon})\leq
T\right)&=\mathbb{P}\left(\{\mathcal{T}_R(u_{\epsilon})>S\}\setminus\{\mathcal{T}_R(u_{\epsilon})>T\}\right)\\
&=\mathbb{P}\left(\mathcal{T}_R(u_{\epsilon})>S\right)
\left(1-\frac{\mathbb{P}\left(\mathcal{T}_R(u_{\epsilon})>T\right)}{\mathbb{P}\left(\mathcal{T}_R(u_{\epsilon})>S\right)}\right).\end{align*}
\end{proof}

\section{Applications to nonlinear optics}\label{s6}
The NLS equation when $d=1$, $\s=1$ and $\lambda=1$ is called the
noisy cubic focusing nonlinear Schr\"odinger equation. It is a
model used in nonlinear optics. Recall that for the above values
of the parameters the solutions are global. The variable $t$
stands for the one dimensional space coordinate and $x$ for the
time. The deterministic equation is such that there exists a
particular class of solutions, which are localized in space (here
time), that propagate at a finite constant velocity and keep the
same shape. These solutions are called solitons or solitary waves.
The functions $$\Psi_{\eta}(t,x)=\sqrt{2}\eta \exp\left(-i\eta^2
t\right)\sech(\eta x),\ \eta>0,$$ form a family of solitons. They
are used in optical fibers as information carriers to transmit the
datum 0 or 1 at high bit rates over long distances. The noise
stands for the noise produced by in-line amplifiers.\vspace{0.3cm}

Let $u_{\epsilon}$ denote the solution with
$u_0(\cdot)=\Psi_1(0,\cdot)$ as initial datum and $\epsilon$ as
noise intensity like in section \ref{s3} and $u_{\epsilon}^n$
denote the solution with null initial datum and the same noise
intensity. The square of the momentum of $u_0$ is
$4$.\vspace{0.3cm}

At a particular coordinate $T$ of the fiber, when a window
$[-l,l]$ is given, the square of the $\xLtwo(-l,l)-$norm, or
measured square of the momentum, is recorded. It is close to the
momentum in the deterministic case for sufficiently high $l$ since
the wave is localized. A decision criterium is to accept that we
have 1 if the measured square of the momentum is above a certain
threshold and 0 otherwise. We set a threshold of the form
$4(1-\gamma)$, where $\gamma$ is a real number in
$[0,1]$.\vspace{0.3cm}

As the soliton is progressively distorted by the noise, it is
possible either to wrongly decide that the source has emitted a 1,
or to wrongly discard a 1. The two error probabilities consist of
$$\mathbb{P}_{\epsilon}^{|0}=\mathbb{P}\left(\int_{-l}^l|u_{\epsilon}^n(T,x)|^2dx\geq 4(1-\gamma)\right)$$
and
$$\mathbb{P}_{\epsilon}^{|1}=\mathbb{P}\left(\int_{-l}^l|u_{\epsilon}(T,x)|^2dx<4(1-\gamma)\right).$$

In the following we make the assumption that $\xker\Phi^*=\{0\}$.
Indeed, from the arguments used in the proof of Proposition
\ref{p3}, it is needed for controllability issues to guaranty that
the infima are not taken over empty sets. Also $T$ is fixed,
$\gamma_0\in \left(0,\frac{1}{2}\right)$ is fixed and the size $l$
of the window is such that
$$\int_{-l}^l|u_d(T,x)|^2dx\wedge\int_{-l}^l|\Psi_1(0,x)|^2dx>4\left(1-\frac{\gamma_0}{2}\right).$$

\begin{prpstn}\label{p6}
For every $\gamma$ in $[\gamma_0,1-\gamma_0]$ besides an at most
countable set of points, the following equivalents for the
probabilities of error hold
\begin{align*}
\log\mathbb{P}_{\epsilon}^{|0}&\sim_{\epsilon\rightarrow0}-\frac{1}{2\epsilon}\inf_{h\in
\xLtwo(0,+\infty;\xLtwo(\R)):\int_{-l}^l|\tilde{\mathbf{S}}(h)(T,x)|^2dx\geq 4(1-\gamma)}\left\{\|h\|_{\xLtwo(0,+\infty;\xLtwo(\R))}^2\right\}\\
\log\mathbb{P}_{\epsilon}^{|1}&\sim_{\epsilon\rightarrow0}-\frac{1}{2\epsilon}\inf_{h\in
\xLtwo(0,+\infty;\xLtwo(\R)):\int_{-l}^l|\mathbf{S}(h)(T,x)|^2dx<4(1-\gamma)}\left\{\|h\|_{\xLtwo(0,+\infty;\xLtwo(\R))}^2\right\}
\end{align*}
where $\tilde{\mathbf{S}}(h)$ is the skeleton associated to the
same control problem as $\mathbf{S}(h)$ but with null initial
datum. Both infima are positive numbers.
\end{prpstn}
\begin{proof} The mapping $\varphi$ from $\mathcal{X}_{\infty}$ into
$\mathbb{R}^+$ such that $\varphi(f)=\int_{-l}^l|f(x)|^2dx$ is
continuous. Therefore, the direct image measures
$\left(\varphi_*\mu^{u_{\epsilon}}\right)_{\epsilon\geq0}$ and
$\left(\varphi_*\mu^{u_{\epsilon}^n}\right)_{\epsilon\geq0}$
satisfy LDP of speed $\epsilon$ and good rate functions
respectively
$$I^T(y)=\frac{1}{2}\inf_{h\in
\xLtwo(0,+\infty;\xLtwo(\R)):\int_{-l}^l|\mathbf{S}(h)(T,x)|^2dx=y}\left\{\|h\|_{\xLtwo(0,+\infty;\xLtwo(\R))}^2\right\}$$
and $J^T$ where $\mathbf{S}$ is replaced by $\tilde{\mathbf{S}}$.
Consequently,
$$\forall i\in\{0,1\},\ -L^i(\gamma)\leq\underline{\lim}_{\epsilon\rightarrow0}\epsilon\log\mathbb{P}_{\epsilon}^i\leq\overline{\lim}_{\epsilon\rightarrow0}\epsilon\log\mathbb{P}_{\epsilon}^i\leq-U^i(\gamma)$$
where
\begin{align*}
&L^0(\gamma)=\inf_{y\in(4(1-\gamma),+\infty)}J^T(y),\ \
U^0(\gamma)=\inf_{y\in[4(1-\gamma),+\infty)}J^T(y),\\
&L^1(\gamma)=\inf_{y\in[0,4(1-\gamma))}I^T(y),\ \ \ \ \ \
U^1(\gamma)=\inf_{y\in[0,4(1-\gamma)]}I^T(y). \end{align*}\\
For every $\delta>0$, $U^0(\gamma)\leq L^0(\gamma)\leq
U^0(\gamma-\delta)$ and $U^1(\gamma)\leq L^1(\gamma)\leq
U^1(\gamma+\delta)$ hold.\\ The function $\gamma\mapsto
U^0(\gamma)$ is positive and decreasing. Also, as
$\xker\Phi^*=\{0\}$, there exists a sequence $(h_n^0)_{n\in\N}$ of
functions of $\xLtwo(0,+\infty;\xLtwo(\R))$ that converges to
$$H^0(t)=i\frac{\xdif}{\xdif t}u^0-\Delta u^0-\lambda|u^0|^{2\sigma}u^0$$
where $$u^0(t)=\indic_{t\leq T}\frac{t}{T}\Psi_1(0,\cdot)$$ and by
the continuity proved in section \ref{s51}
$\left(\varphi\circ\mathbf{S}(h_n^0)\right)_{n\in\N}$ converges to
$\varphi\circ\mathbf{S}(H^0)$ which is such that
$\varphi\circ\mathbf{S}(H^0)>4\left(1-\frac{\gamma_0}{2}\right)>4(1-\gamma_0)$.
Consequently, $h_n^0$ belongs to the minimizing set for $n$ large
enough. Thus, $U^0(\gamma_0)<+\infty$ follows. Consequently, the
function $\gamma\mapsto U^0(\gamma)$ possesses an at most
countable set of points of discontinuity.\\
Similarly, the function $\gamma\mapsto U^1(\gamma)$ is a bounded
increasing function. Also, if $(h_n^1)_{n\in\N}$ and $H^1(t)$ are
defined as previously replacing $u^0(t)$ by
$$u^1(t)=\indic_{t\leq T}\left(1-\left(1-\sqrt{\frac{\gamma_0}{2}}\right)\frac{t}{T}\right)\Psi_1(0,\cdot),$$
the sequence $(\varphi\circ\mathbf{S}(h_n^1))_{n\in\N}$ converges
to $\varphi\circ\mathbf{S}(H^1)\leq
2\gamma_0=4\left(1-\left(1-\frac{\gamma_0}{2}\right)\right)$.
Thus, for $n$ large enough $h_n^1$ belongs to the minimizing set.
Consequently, the function $\gamma\mapsto U^1(\gamma)$ has an at
most countable set of points of discontinuity. Thus, for a well
chosen $\gamma$, letting $\delta$ converge to zero, we obtain for
$i\in\{0,1\}$ that $L^i(\gamma)=U^i(\gamma)$ and the

equivalents follow.\\
From the arguments used in the proof of Proposition \ref{p4},
$\tilde{\mathbf{S}}$ is a continuous mapping from
$\xLtwo(0,+\infty;\xHone(\R))$ into $\mathcal{X}_{\infty}$. Since
$\varphi$ is continuous, if $(H_n)_{n\in\N}$ is a sequence of
functions converging to zero in $\xLtwo(0,+\infty;\xHone(\R))$
then $\left(\varphi\circ\tilde{\mathbf{S}}(H_n)\right)_{n\in\N}$
converges to $\varphi\circ\tilde{\mathbf{S}}(0)=0$. Proposition
\ref{p4} also gives that
$\left(\varphi\circ\mathbf{S}(H_n)\right)_{n\in\N}$ converges to
$\varphi\circ\mathbf{S}(0)$ which satisfies
$\varphi\circ\mathbf{S}(0)>4\left(1-\frac{\gamma_0}{2}\right)$.
The conclusion follows.
\end{proof}

In reference \cite{FKLT} the authors explain that, for the second
error probability, two processes are mainly responsible for the
deviations of the measured square of the momentum from its
expected value: the fluctuation of the soliton power
$$\frac{M(u_{\epsilon}(T))^2}{4}$$ and a shift of the soliton
position, also called center of mass, timing jitter or fluctuation
in timing since time and space variables have been exchanged,
characterized by
$$\frac{\int_{-\infty}^{+\infty}x|u_{\epsilon}(T,x)|^2dx}{M(u_{\epsilon}(T))^2}.$$
Recall that $M$ stands for the momentum or $\xLtwo(\R)-$norm. The

authors give an asymptotic expression of the probability density
function of the joint law of the two above random variables. In
the case of the first error probability, they explain that the
optimal way to create a large signal is to grow a soliton, they
thus give the marginal probability density function of the square
of the momentum. \vspace{0.3cm}

In the two following sections we concentrate on the square of the
momentum, we take $l=+\infty$ as if the window were not bounded.
Somehow, if we forget the coefficient, we concentrate on the tails
of the marginal law of the soliton power when $\epsilon$ converges
to zero. We recall, it has been pointed out in the introduction,
that the momentum is no longer preserved in the stochastic case
and is such that its expected value increases. We then study the
tails of the law of the shift of the soliton position when
$\epsilon$ converges to zero. Remark that we drop the
renormalization in the shift of the soliton position so as to
obtain a probability measure in the integral and write
$$Y_{\epsilon}=\int_{-\infty}^{+\infty}x|u_{\epsilon}(T,x)|^2dx.$$
We finally present a result for the general case, with no
limitation on $\sigma$ and $d$, where blow-up may occur.

\subsection{Upper bounds}\label{s61}
The norm of the linear continuous operator $\Phi$ of $\xLtwo(\R)$
is thereafter denoted by
 $\|\Phi\|_c$.
\begin{prpstn}\label{p7}
For every positive $T$, $\gamma$ in $[0,1]$, and every operator
$\Phi$ in\\ $\mathcal{L}_2(\xLtwo(\R),\xHone(\R))$, the
inequalities
$$\overline{\lim}_{\epsilon\rightarrow0}\epsilon\log\mathbb{P}_{\epsilon}^{|0}\leq-\frac{1-\gamma}{2T\|\Phi\|_c^2}$$

and
$$\overline{\lim}_{\epsilon\rightarrow0}\epsilon\log\mathbb{P}_{\epsilon}^{|1}\leq-\frac{1+\gamma}{T\|\Phi\|_c^2}\left(\sqrt{1+\left(\frac{\gamma}{1+\gamma}\right)^2}-1\right).$$
hold.\end{prpstn}
\begin{proof} Multiplying by $-i\overline{u}$ the equation
$$i\frac{\xdif}{\xdif t}u-\Delta u-\lambda|u|^{2\sigma}u=\Phi h,$$
integrating over time and taking the real part gives that
$$\|u(T)\|_{\xLtwo(\R)}^2-\|u_0\|_{\xLtwo(\R)}^2=2\Re\left(-i\int_0^T\int_{\R}\Phi h \overline{u}\ dxdt\right).$$

\textbf{First bound: }The boundary conditions
$\|u(T)\|_{\xLtwo(\mathbb{R})}^2\geq 4(1-\gamma)$ and $u_0=0$
along with Cauchy Schwarz inequality imply both that
$$4(1-\gamma)\leq 2\|\Phi\|_c\|h\|_{\xLtwo(0,T;\xLtwo(\R))}\|u\|_{\xLtwo(0,T;\xLtwo(\R))}
,$$ and that
\begin{align*}
\int_0^T\|u(t)\|_{\xLtwo(\R)}^2dt&=2\int_0^T\Re\left(-i\int_0^t\Phi
h \overline{u}\ dxds\right)dt\\
&\leq 2
T\|\Phi\|_c\|h\|_{\xLtwo(0,T;\xLtwo(\R))}\|u\|_{\xLtwo(0,T;\xLtwo(\R))},
\end{align*}
thus,
$$\|h\|_{\xLtwo(0,+\infty;\xLtwo(\R))}^2\geq\frac{1-\gamma}{T\|\Phi\|_c^2}.$$
\textbf{Second bound: }The new boundary conditions
$\|u(T)\|_{\xLtwo(\R)}^2<4(1-\gamma)$ and
$\|u_0\|_{\xLtwo(\R)}^2=4$ give both that along with Cauchy
Schwarz inequality
$$4\gamma<2\|\Phi\|_c\|h\|_{\xLtwo(0,+\infty;\xLtwo(\R))}\|u\|_{\xLtwo(0,T;\xLtwo(\R))}$$
and also along with Cauchy Schwarz and integration over time
$$\|u\|_{\xLtwo(0,T;\xLtwo(\R))}^2-4T\leq2\|\Phi\|_c\|h\|_{\xLtwo(0,+\infty;\xLtwo(\R))}\|u\|_{\xLtwo(0,T;\xLtwo(\R))}.$$
Consequently, it follows that
$$\|u\|_{\xLtwo(0,T;\xLtwo(\R))}\leq
T\|\Phi\|_c\|h\|_{\xLtwo(0,T;\xLtwo(\R))}\left(1+\sqrt{1+\frac{4}{T\|\Phi\|_c^2\|h\|_{\xLtwo(0,T;\xLtwo(\R))}^2}}\right).$$
Thus, we obtain
$$\frac{2\gamma}{T\|\Phi\|_c^2}<\|h\|_{\xLtwo(0,+\infty;\xLtwo(\R))}^2\left(1+\sqrt{1+\frac{4}{T\|\Phi\|_c^2\|h\|_{\xLtwo(0,T;\xLtwo(\R))}^2}}\right)$$
and
$$\|h\|_{\xLtwo(0,T;\xLtwo(\R))}^2>\frac{2(1+\gamma)}{T\|\Phi\|_c^2}\left(\sqrt{1+\left(\frac{\gamma}{1+\gamma}\right)^2}-1\right).$$

The upper bound follows. \end{proof}

\subsection{Lower bounds} \label{s63}
We prove the following lower bounds.
\begin{prpstn}\label{p8}
For every positive $T$, $\gamma$ in $[0,1]$, and every operator
$\Phi$ in\\
$\mathcal{L}_2(\xLtwo(\R),\xHone(\R))$ acting as the identity map
on $\spn\{\frac{1}{\ch(a x)},x\frac{\sh}{\ch^2}(a x);a\in\R\}$,
the inequalities
$$\underline{\lim}_{\epsilon\rightarrow0}\epsilon\log\mathbb{P}_{\epsilon}^{|0}\geq-\frac{2(1-\gamma)(12+\pi^2)}{9T}$$
and
$$\underline{\lim}_{\epsilon\rightarrow0}\epsilon\log\mathbb{P}_{\epsilon}^{|1}\geq-\frac{2(2-\gamma-2\sqrt{1-\gamma})(12+\pi^2)}{9T}.$$
hold.
\end{prpstn}
\begin{proof}
Let the constant $\eta$ in the parametrized family of solitons
depend on $t$ and set
\begin{equation}
\label{e6}u(t,x)=\Psi_{S}(t,x)=\sqrt{2}\eta(t)
\exp\left(-i\int_0^t\eta^2(s)ds\right)\sech(\eta(t) x)
\end{equation}
then, from the assumption on $\Phi$, the function

$$h_{S}(t,x)=i\frac{\eta'(t)}{\eta(t)}\Psi_{S}(t,x)-i\sqrt{2}\eta'(t)\eta(t)x\exp\left(-i\int_0^t\eta^2(s)ds\right)\frac{\sh\left(\eta(t)x\right)}{\ch^2\left(\eta(t)x\right)}$$
is such that $\mathbf{S}(h_{S})$ and $\tilde{\mathbf{S}}(h_{S})$
are the solutions of the control problems. Also, as $u_0$ belongs
to $\xHtwo(\R)$, $\mathbf{S}(h_{S})$ and
$\tilde{\mathbf{S}}(h_{S})$ are functions of
$\xCn([0,T];\xHtwo(\R))\cap\xCn^1([0,T];\xLtwo(\R))$, consequently
$t\rightarrow\eta(t)=\frac{1}{4}\|\Psi_S(t,\cdot)\|_{\xLtwo(\R)}^2$
is necessarily a function in $\xCn^1([0,T])$.\\
For the first error probability, the lower bound follows from the

fact that the infimum is smaller than the infimum on the smallest
set of parametrized $h_S$ and the computation of the
$\xLtwo(0,T;\xLtwo(\R))$ norm of $h_S$ which gives that
$$\inf_{\eta\in\
\xCn^1([0,T]):\eta(0)=0,\|\tilde{\mathbf{S}}(h_S)(T,\cdot)\|_{\xLtwo(\R)}^2dx\geq4(1-\gamma)}\left\{\|h_S\|_{\xLtwo(0,+\infty;\xLtwo(\R))}^2\right\}$$is
equal to$$\inf_{\eta\in\
\xCn^1([0,T]),b.c.}\int_0^TF_{S}(\eta(t),\eta'(t))dt,$$ where the
Lagrangian $F_S$ is
$$F_S(z,p)=\frac{1}{9}(12+\pi^2)\frac{p^2}{z},$$
and b.c. stands for the boundary conditions $\eta(0)=0$ and
$\eta(T)\geq1-\gamma$. Indeed, since $\tilde{\mathbf{S}}(h)(T)$ is
a function of $(h(t))_{t\in[0,T]}$, the infimum could be taken on
functions set to zero almost everywhere after $T$, thus
$\|h\|_{\xLtwo(0,+\infty;\xLtwo(\R))}^2$ in the left hand side
could be replaced by $\|h\|_{\xLtwo(0,T;\xLtwo(\R))}^2$. A scaling
argument gives that the terminal boundary condition is
necessarily saturated.\\
Similarly, for the second error probability, $\tilde{\mathbf{S}}$
is replaced by $\mathbf{S}$ and b.c. is $\eta(0)=1$ and
$\eta(T)=1-\gamma$.\\
The usual results of the indirect method do not apply to the
problem of the calculus of variations, nonetheless solutions of
the boundary value problem associated to the Euler-Lagrange
equation
$$2\frac{\eta''}{\eta}=\left(\frac{\eta'}{\eta}\right)^2$$ provide
upper bounds when we compute the integral of the Lagrangian. If we
suppose that $\eta$ is in $\xCn^{3}([0,T])$ and that it is

positive on $(0,T)$, we obtain by derivation of the ODE that on
$(0,T)$,
$$\eta'''=0.$$ Also, looking for solutions of the form
$at^2+bt+c$, we obtain that necessarily $b^2=4ac$. Thus
$\xCn^{3}([0,T])$ positive solutions are necessarily of the form
$a\left(t-\frac{b}{2a}\right)^2$. >From the boundary conditions, we
obtain that for the first error probability the function defined
by
$$\eta^0(t)=(1-\gamma)\left(\frac{t}{T}\right)^2$$
is a solution of the boundary value problem. For the second error
probability, the boundary conditions imply that the two following

functions defined by
$$\eta^{1,1}(t)=\left(2-\gamma+2\sqrt{1-\gamma}\right)\left(\frac{t}{T}\right)^2+2\left(-1-\sqrt{1-\gamma}\right)\frac{t}{T}+1$$
and
$$\eta^{1,2}(t)=\left(2-\gamma-2\sqrt{1-\gamma}\right)\left(\frac{t}{T}\right)^2+2\left(-1+\sqrt{1-\gamma}\right)\frac{t}{T}+1$$

are solutions of the boundary value problem. The second function
gives the smallest value when we compute the integral of the
Lagrangian.
\end{proof}
Remark that, in the case of the first error probability, both
upper and lower bounds in Proposition \ref{p7} and Proposition
\ref{p8} are increasing functions of $\gamma$. Similarly, in the
case of the second error probability, the bounds are decreasing
functions of $\gamma$. This could be interpreted as the higher is
the threshold, the more energy is needed to form a signal which
momentum gets above the threshold at the coordinate $T$ and
conversely in the case of a soliton as initial datum.
\begin{rmrk}
When there is no particular tradeoff between the two errors, the
overall risk of error in transmission, due to noise, can be taken
as
$$R_{\epsilon}=\mathbb{P}_{\epsilon}^{|0}\vee\mathbb{P}_{\epsilon}^{|1}.$$
Choosing $\gamma=\frac{5}{7}$ allows to minimize the maximum of
the two upper bounds of Proposition \ref{p7}, for the associated

threshold we get
$$\overline{\lim}_{\epsilon\rightarrow0}\epsilon\log
R_{\epsilon}\leq-\frac{1}{7T\|\Phi\|_c^2}.$$ Similarly, choosing
$\gamma=\frac{3}{4}$ minimizes the maximum of the two lower bounds
of Proposition \ref{p8}, for the associated threshold we get
$$-\frac{12+\pi^2}{18T}\leq\underline{\lim}_{\epsilon\rightarrow0}\epsilon\log
R_{\epsilon}.$$
\end{rmrk}
The two values of $\gamma$ are very close and correspond to
thresholds taken as $1.14$ and $1$ respectively, i.e. to $22\%$
and $25\%$ of the momentum of the initial datum.\vspace{0.3cm}

Remark also that the bounds for the error probabilities are of the
right order. Indeed, from the probability density function given
for the first error probability in \cite{FKLT}, we are expecting,
when the noise is the ideal white noise and thus $\|\Phi\|_c=1$,
that
$\lim_{\epsilon\rightarrow0}\epsilon\log\mathbb{P}_{\epsilon}^{|0}=-c\frac{1-\gamma}{T}$
with a positive constant $c$. The joint probability density
function obtained in \cite{FKLT}, when the initial datum is a
soliton, also allows to obtain asymptotics of the probability of
tail events of the soliton power and thus of the second error
probability.\vspace{0.3cm}

Remark finally that it is natural to obtain that the opposite of
the error probabilities are decreasing functions of $T$. Indeed,
the higher is $T$, the less energy is needed to form a signal
which momentum gets above a fixed threshold at the coordinate $T$.
Replacing above by under, we obtain the same result in the case of
a soliton as initial datum. Consequently, the higher is $T$ the
higher the error probabilities get.\vspace{0.3cm}

We give in the next section some further considerations on
problems of the calculus of variations along with the results of
some numerical computations.

\subsection{Remarks on the problem of
the calculus of variations}\label{s64} The most efficient
parametrization we obtained has been presented above. Many
parametrizations are at hand and we may expect one that allows to
obtain a good approximation of the rate of convergence to zero of
the error probabilities. This would allow to have an idea of the
most likely trajectories leading to an error in transmission. The
presence of an operator $\Phi$ before $h$ in the control problem
is a limitation for the computation. We will thereafter consider
that $\Phi$ acts as the identity map on a sufficiently large
linear space of $\xLtwo(\R)$.\vspace{0.3cm}

\textbf{A parametrization as a function of $t$ and $x$.} If we
look for solutions of the form
$u_0\left(1-\frac{t}{T}\right)+u_T\frac{t}{T}$ and suppose that
$u_T$ belong to $\xHn^{3}(\R^d)$, the infimum over the functions
$u_T$ gives birth to a Euler-Lagrange equation which is a fourth
order nonlinear PDE. A more reasonable approach for computation is
to find an upper bound as an infimum functions parametrized by

certain paths on $\R$.\vspace{0.3cm}

\textbf{The parametrization by the amplitude.} We consider paths
of the form $u(t,x)=f(t)u_0(x)$ for $f\in\xHone(0,T)$. In that
case $\xHone(0,T)$ is a space of real valued functions. We recall
that, from the Sobolev inequalities, $\xHone(0,T)$ is continuously
embedded into the space of $\frac{1}{2}-$H\"older functions
$\xCn^{\frac{1}{2}}([0,T])$ and thus into every $\xLn^{p}(0,T)$.
Remark that the function $h_a(t)=if'(t)u_0-f(t)\Delta
u_0-f(t)^3|u_0|^2u_0$ is the associated control. Minimizing over
the norm of these controls leads to a problem of the calculus of
variations. We obtain that
$$L^{(t,+\infty)}\leq\frac{1}{2}\inf_{f\in\mathcal{A}}\int_0^TF_a(f(t),f'(t))dt$$
where the Lagrangian $F_a$ is defined by
$$F_a(z,p)=p^2\|u_0\|_{\xLtwo(\R)}^2+z^2\|\Delta
u_0\|_{\xLtwo(\R)}^2+2 z^4\left(\Delta
u_0,|u_0|^2u_0\right)_{\xLtwo(\R)}+z^6\||u_0|^2u_0\|_{\xLtwo(\R)}^2$$
which, with the particular value of $u_0$, becomes
$$F_a(z,p)=4p^2+\frac{28}{15}z^2-\frac{32}{5}z^4+\frac{128}{15}z^6,$$
and $\mathcal{A}$ is the admissible set
$\{f\in\xHone(0,T):f(0)=0,\ f(T)=\sqrt{1-\gamma}\}$ for the first
error probability and $\{f\in\xHone(0,T):f(0)=1,\
f(T)=\sqrt{1-\gamma}\}$ for the second error probability.\\
Indeed, as for every $\alpha\geq1$,
$\frac{28}{15}(\alpha^2-1)z^2-\frac{32}{5}(\alpha^4-1)z^4+\frac{128}{15}(\alpha^6-1)z^6\geq0$,
we have $F_a(\alpha z,\alpha p)\geq F_a(z,p)$, thus the terminal
boundary condition corresponding to the first error probability is
necessarily saturated. The same holds for the terminal value
corresponding to the second error probability changing $z$ to
$z-1$.\vspace{0.3cm}

Also, standard calculation gives that for every $(z,p)$ in $\R^2$,
$F_a(z,p)\geq 4p^2+\frac{2}{3}z^2$. Recall that as $F_a$ is
non-negative, convex in the $p$ variable, $F_a$ and
$\frac{\partial}{\partial p}F_a$ are continuous in the $(z,p)$
variables, the integral is a weakly sequentially lower
semicontinuous function of $f$. In addition, as in both cases
$\mathcal{A}$ is nonempty and $F_a$ satisfies the coercivity
condition $F_a(z,p)\geq4p^2$ for every $(z,p)$ in $\R^2$, there
exists at least one minimizer, i.e. a function that solves the

problem of the calculus of variations.\vspace{0.3cm}

It is finally also possible to revisit slightly the proof of
Proposition $4$ of section $8.2.3$ of reference \cite{ELC} and to
check that, in our particular case, though we do not have the
expected growth conditions on the Lagrangian, Lebesgue's dominated
convergence theorem could be applied and that the result of the
proposition still holds. Consequently, any minimizer is a weak
solution of the Euler-Lagrange equation. If we apply the chain
rule supposing that $f$ is a function of $\xHtwo(0,T)$, the
Euler-Lagrange equation becomes the following nonlinear ODE
$$15f''-7f+48f^3-96f^5=0.$$
By numerical computations of the integral and of the boundary
value problems for the first and second error probabilities, we
obtain in figure $1$ the curves labeled "amplitude". These bounds
are compared to the bounds of section \ref{s63} where the curves
are labeled "soliton parameter" and to
the upper bounds of section \ref{s61}.\\
\begin{figure}[H]
\centering
\includegraphics[height=0.44\hsize]{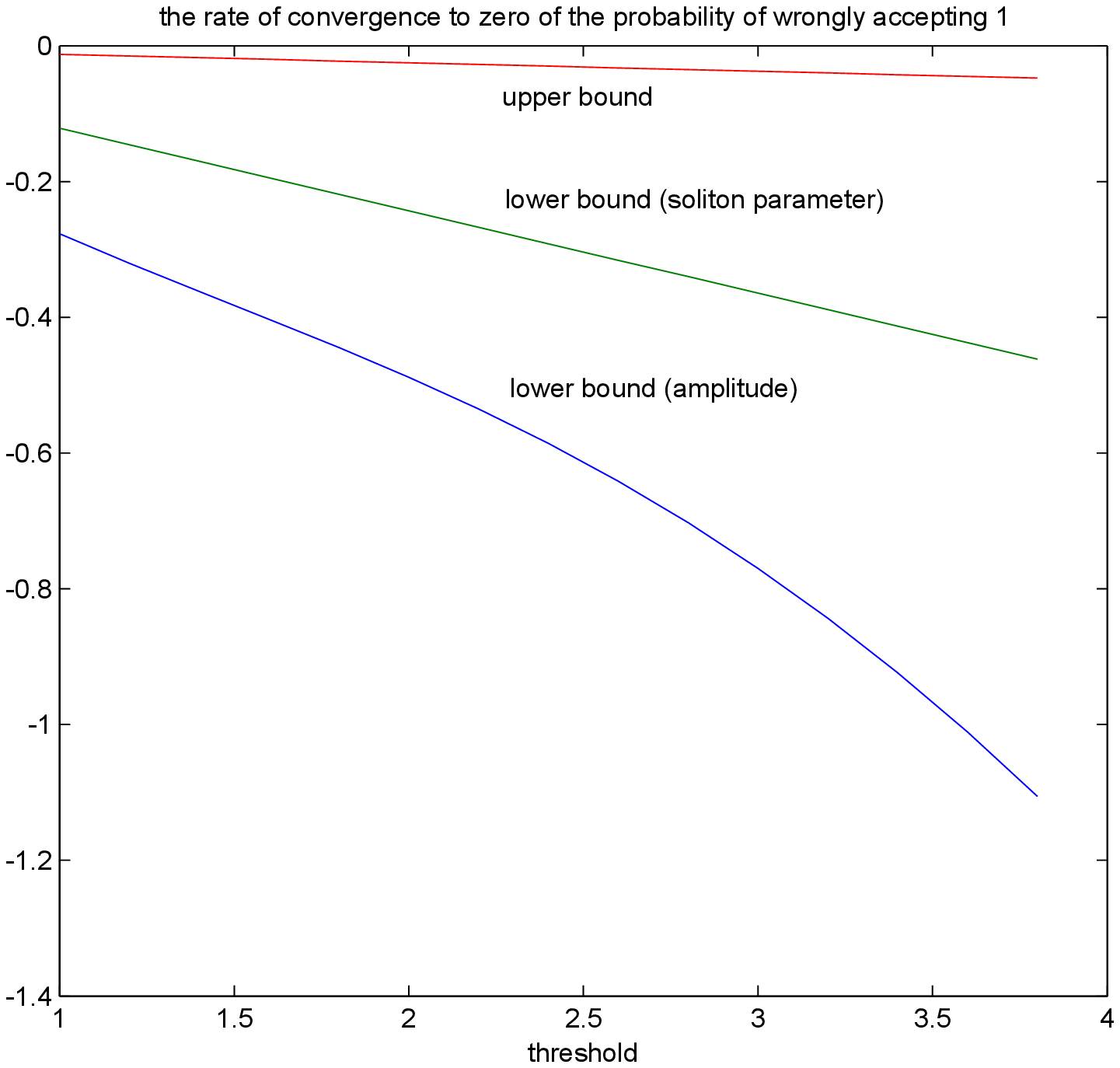}
\qquad
\includegraphics[height=0.44\hsize]{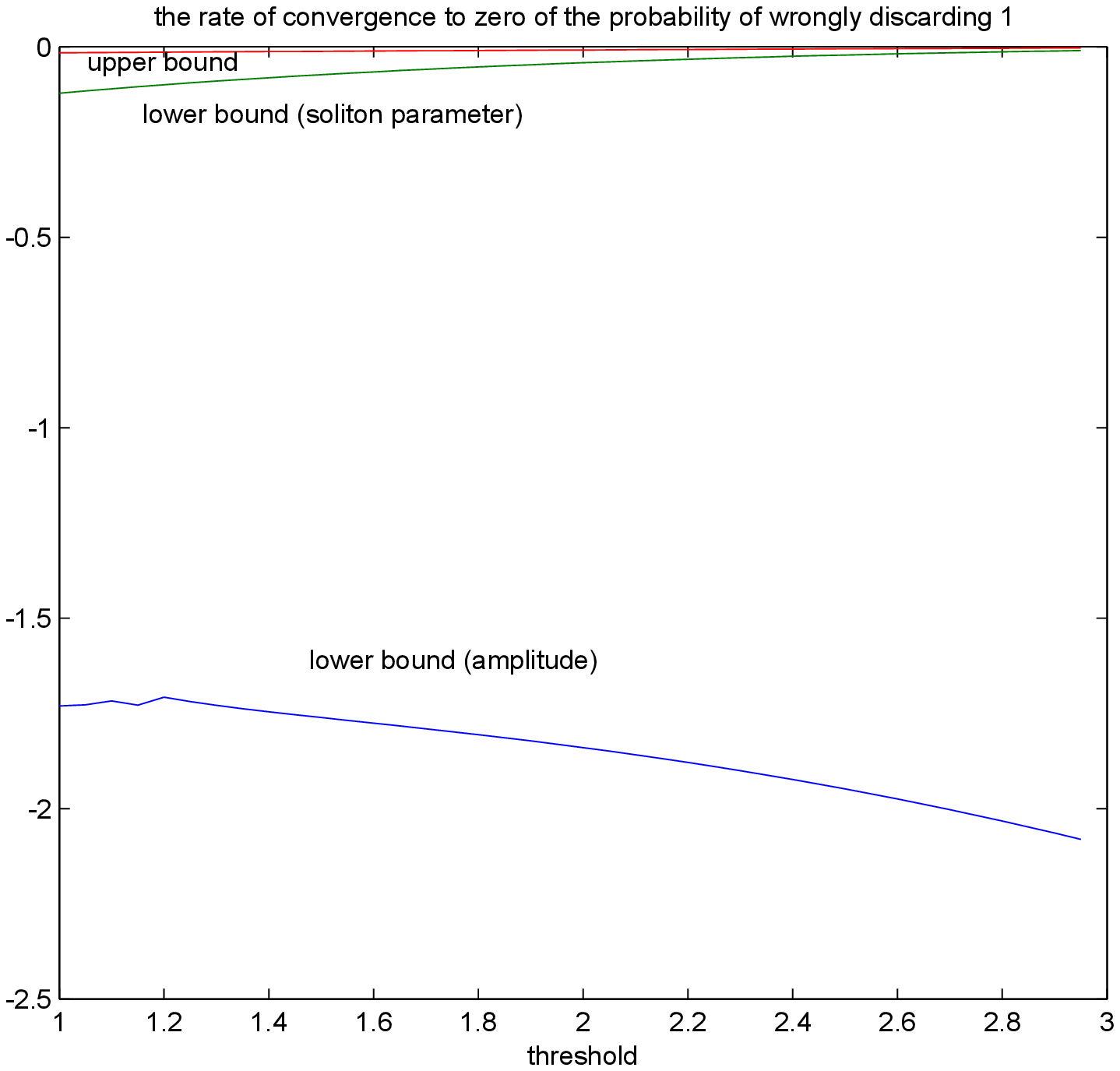}
\caption{Bounds as a function of the threshold $4(1-\gamma)$ with
$T=10$ and $\|\Phi\|_c=1$.}
\end{figure}
As we may see in figure $1$, this parametrization is less
efficient to obtain lower bounds. Moreover, the case of interest
in fibers is when $T$ is large and a graph of $T$ times the
integral shows that this parametrization gives a lower bound of

order less than $\frac{1}{T}$. It is not satisfactory, nonetheless
it is presented in the paper because it could be used with
arbitrary $u_0$, $\sigma$ and $d$, to obtain some refinement of
the result of Proposition \ref{p3}. Remark also that in this case
a nice parametrization could give an idea of the most likely
trajectories of the solutions that blow up after time
$T$.\vspace{0.3cm}

These computations show that a parametrization of the amplitude of
the soliton is not enough to approximate the most likely
trajectories leading to an error in transmissions. The phase and
several other shape parameters have to be introduced. Remark that

we also tried a parametrization replacing in the exponential in
the identity (\ref{e6}) the term $\int_0^t\eta^2(s)ds$ by
$\eta^2(t)t$ but that it gave less interesting results.
\vspace{0.3cm}

\textbf{A more complete parametrization.} In reference \cite{FKLT}
the authors parametrize the signal in the following way
$$u(t,x)=\Psi_p(t,x)=\sqrt{2}\eta(t) \exp\left(i\beta(t)
x+i\alpha(t)+i\tau(t)\right)\left[\sech(\eta(t)(x-y(t)))+v(x)\right],$$
where $\tau$ is the "internal time", it satisfies

$\tau'(t)=\eta^2(t)$, and $v$ the "continuous spectrum on the
background of the soliton". This parametrization is used to obtain
a good approximation of the probability density function of the
joint law of the soliton power and of the shift of the soliton
position. We saw that introducing a field $v$ strongly complicates
the calculus of variations. An interested reader could refer to
the article for a physical justification that $v$ is indeed
negligible in first approximation. The authors also neglect the
$\alpha$ variable, we keep it. The associated control is the

function
\begin{align*}
&h_p(t,x)=i\frac{\eta'(t)}{\eta(t)}\Psi_p(t,x)-(\beta'(t)x+\alpha'(t))\Psi_p(t,x)\\
&+\sqrt{2}i[-\eta'(t)(x-y(t))+\eta(t)y'(t)]\exp\left(i\beta(t)x+i\alpha(t)+i\tau(t)\right)\eta(t)\frac{\sh}{\ch^2}(\eta(t)(x-y(t))).
\end{align*}
The augmented Lagrangian of the new problem of the calculus of
variations is given by
$$F_P(Z,P)=\frac{12+\pi^2}{9}\frac{p_1^2}{z_1}+\frac{4}{3}z_1^3p_4^2+4p_2^2z_1+4p_3^2z_1z_4^2+\frac{\pi^2p_3^2}{3z_1}+8p_2p_3z_1z_4,$$
where $Z=(z_1,z_2,z_3,z_4)=(\eta,\alpha,\beta,y)$ and
$P=(p_1,p_2,p_3,p_4)$ is its derivative. A scaling argument still
indicates that the terminal boundary conditions are saturated,
thus $\eta(0)=0$ and $\eta(T)=1-\gamma$ in the case of the first
error probability and $\eta(0)=1,\ \alpha(0)=\beta(0)=y(0)=0$ and
$\eta(T)=1-\gamma$ for the second error probability. Still, the
usual results of the indirect method do not apply and we are not
assured that a minimizer is a solution of the Euler-Lagrange nor

of the contrary. Nonetheless, the solutions of the Euler-Lagrange
equation may be good guesses and provide upper bounds of the
problem of the calculus of variations and thus lower bounds of the
rate of exponential convergence to zero of the error
probabilities. It consists of a system of coupled nonlinear ODEs.
The natural boundary conditions at the free end leads to
$\alpha'=\beta'=0$ and the system simplifies in
$$\left\{\begin{array}{c}
    \frac{12+\pi^2}{9}\left(\frac{2\eta''}{\eta}-\left(\frac{\eta'}{\eta}\right)^2\right)=4\eta^2(y')^2\\
    \eta'\eta^2y'+\frac{\eta^3y''}{3}=0\\
\end{array}\right.$$with the extra condition $y'(T)=0$ resulting from the natural boundary conditions. It is associated to the simpler Lagrangian
$$F_{PS}(z_1,z_4,p_1,p_4)=\frac{12+\pi^2}{9}\frac{p_1^2}{z_1}+\frac{4}{3}z_1^3p_4^2.$$
The singularity did not allow to treat numerically the case of
null initial datum and led to an identically null solution for
$y'$ and thus to the same solution for $\eta$ and the same lower
bound as in section \ref{s63} in the case of a soliton as initial
datum.\vspace{0.3cm}

The parametrization is probably very well adapted to the case of a
bounded window, particularly for the second error probability. The
variable $y$ is a parametrization of the shift of the soliton
position, it is also likely to be relevant for the computation of
a lower bound by the calculus of variations in the next section.

\subsection{The tails of the law of the shift of
the soliton position when $\epsilon$ converges to zero} We
consider the mapping $Y$ from $\Sigma^{\frac{1}{2}}$ into $\R$
defined by $$Y(f)=\int_{\R}x|f(x)|^2dx.$$ We suppose that $\Phi$
belongs to $\mathcal{L}_2(\L2,\Sigma)$ and that $\xker
\Phi^*=\{0\}$. The following lemma allows to define the family of
measures $\left(\mu^{Y_{\epsilon}}=(Y\circ
eval_T\circ\mathcal{G})_*\mu^{Z_{\epsilon}}\right)_{\epsilon>0}$
which correspond to the laws of the shift of the position of the
soliton for each solution $u_{\epsilon}$ with a soliton as initial
datum and for a noise of intensity $\epsilon$. Also, since the
shift is zero for the deterministic solution $u_d$, the lemma
gives that the sequence of measures converges weakly to the Dirac
mass on $0$. Proposition \ref{p9} characterizes the convergence to
zero of the tails.
\begin{lmm}
For every positive $T$, the mapping $Y\circ\ eval_T\circ\
\mathcal{G}$ from\\
$X^{(T,4);\Sigma}=\xCn([0,T];\Sigma)\cap\xL^{r(4)}\left(0,T;\xW^{1,4}(\R)\right)$
into $\R$ is continuous.
\end{lmm}
\begin{proof}
Let $Z$ and $Z'$ belong to $X^{(T,4);\Sigma}$, the triangular
inequality along with H\"older's inequality  allow to compute the
following sequence of inequalities
$$\left|\int_{\R}x\left(|\mathcal{G}(Z)(T,x)|^2-|\mathcal{G}(Z')(T,x)|^2\right)dx\right|$$
$$\leq\int_{\R}|x|(|\mathcal{G}(Z)(T,x)|+|\mathcal{G}(Z)(T,x)|)|(|\mathcal{G}(Z)(T,x)|-|\mathcal{G}(Z)(T,x)|)|dx$$
$$\leq\|\mathcal{G}(Z)(T)-\mathcal{G}(Z')(T)\|_{\xLtwo(\R)}^{\frac{1}{2}}\sqrt{V(|\mathcal{G}(Z)(T)|+|\mathcal{G}(Z')(T)|)},$$
where $V$ is the variance , it is defined in section \ref{s2}.\\
After some calculations we obtain that
$$\sqrt{V(|\mathcal{G}(Z)(T)|+|\mathcal{G}(Z')(T)|)}$$ is lower than $$2\sqrt{2}\left(\sqrt{V(v(Z)(T))}+\sqrt{V(v(Z')(T))}+\sqrt{V(Z(T))}+\sqrt{V(Z'(T))}\right).$$
The application of Gronwall's inequality given in the proof of

Proposition $3.5$ of reference \cite{dBD2}, along with the Sobolev
injection of $\xHone(\R)$ into $\xLinfty(\R)$ and the continuity
of $\mathcal{G}$ from $X^{(T,4);\Sigma}$ into
$\xCn([0,T];\xHone(\R))$ give that this last term is bounded when
$Z$ and $Z'$ are in a bounded set. We conclude using the
continuity of $\mathcal{G}.$\end{proof} The fact that $Z$ defines
a $X^{(T,4);\Sigma}$ random variable follows from similar
arguments as those used in the proof of Proposition \ref{p1}. We
still denote by $\mu^Z$ its law and by $\mu^{Z_{\epsilon}}$ the
direct images under the transformation $x\mapsto\sqrt{\epsilon} x$
on $X^{(T,4);\Sigma}$.
\begin{prpstn}\label{p9}
The family of measures $(\mu^{Y_{\epsilon}})$ on $\R$ satisfy a
LDP of speed $\epsilon$ and good rate function
$$I^Y(y)=\frac{1}{2}\inf_{h\in\xLtwo(0,T;\xLtwo(\R)):\int_{\R}x|\mathbf{S}(h)(T,x)|^2dx=y}\left\{\|h\|_{\xLtwo(0,T;\xLtwo(\R))}^2\right\}.$$
Moreover for every nonempty interval $J$ such that $0\notin J$,
$0<\inf_{y\in J}I^Y(y)<\infty,$\\
and for every positive $R$ besides an at most countable set of
points, the following equivalents hold
\begin{center}
\begin{align*}
\log \mathbb{P}(Y_{\epsilon}\geq
R)&\sim_{\epsilon\rightarrow0}-\frac{1}{2\epsilon}\inf_{h\in
\xLtwo(0,T;\xLtwo(\R)):\int_{\R}x|\mathbf{S}(h)(T,x)|^2dx\geq
R}\left\{\|h\|_{\xLtwo(0,T;\xLtwo(\R))}^2\right\}\\
\log \mathbb{P}(Y_{\epsilon}\leq
-R)&\sim_{\epsilon\rightarrow0}-\frac{1}{2\epsilon}\inf_{h\in
\xLtwo(0,T;\xLtwo(\R)):\int_{\R}x|\mathbf{S}(h)(T,x)|^2dx\leq
-R}\left\{\|h\|_{\xLtwo(0,T;\xLtwo(\R))}^2\right\}.
\end{align*}
\end{center}
\end{prpstn}
\begin{proof}The LDP for the family $(\mu^{Z_{\epsilon}})_{\epsilon>0}$, which
are centered gaussian measures on a real Banach space, the fact
that their RKHS is $\xim \mathcal{L}$ with the norm of the image
structure and Varadhan's contraction principle give that the
family $\left(\mu^{Y_{\epsilon}}\right)_{\epsilon>0}$ satisfy a
LDP of speed $\epsilon$ and good rate function which is the rate
function of the proposition
$$I^Y(y)=\inf_{z\in X^{(T,4);\Sigma}:Y\circ\ eval_T\circ\
\mathcal{G}(z)=y}\left\{\inf_{h\in\xLtwo(0,T;\xLtwo(\R)):
\mathcal{L}(h)=z}\left\{\frac{1}{2}\|h\|_{\xLtwo(0,T;\xLtwo(\R))}^2\right\}\right\}.$$
The fact that $\inf_{y\in J}I^Y(y)<\infty$ follows from the
assumption $\xker\Phi^*=\{0\}$ and that for every real number $a$,
a solution of the form $u(t,x)=(1+atx)u_0$ satisfies
$Y(u(T))=\frac{aT\pi^2}{3}$, i.e. there exists controls such that the solution reaches any interval at "time" $T$.\\
The positivity follows from the continuity of the mapping $Y\circ
eval_T\circ \mathcal{G}\circ \Lambda\circ \Phi$ from
$\xLtwo(0,T;\xLtwo(\R))$ into $\R$. Indeed $\Phi$ is continuous
from $\xLtwo(0,T;\xLtwo(\R))$ into $\xLone(0,T;\Sigma)$ and
$\Lambda$ is continuous from $\xLone(0,T;\Sigma)$ into
$X^{(T,4);\Sigma}$. This last statement follows from the
Strichartz estimates, the facts that the group is a unitary group
on $\Sigma$ and that the integral is a linear continuous
mapping from $\xLone(0,T;\Sigma)$ into $\Sigma$.\\
The proofs of the final statements are the same as those of
Proposition \ref{p6}.
\end{proof}

Remark that the results are true in $\R^d$ with an arbitrary order
of subcritical nonlinearity.
\subsection{The transmission in the general case where blow-up
may occur} We are interested by estimates of the probabilities
$$\mathbb{P}\left(\int_{-l}^l|u_{\epsilon}^n(T,x)|^2dx\geq4-\gamma,\mathcal{T}(u_{\epsilon}^n)>T\right)=\mathbb{P}_{\epsilon}^{|0,\mathcal{E}}$$
and of
$$\mathbb{P}\left(\int_{-l}^l|u_{\epsilon}(T,x)|^2dx<4-\gamma,\mathcal{T}(u_{\epsilon})>T\right)=\mathbb{P}_{\epsilon}^{|1,\mathcal{E}}.$$
\begin{prpstn}\label{p10}
The following inequalities for the two error probabilities hold,
$$-L^{|0,\mathcal{E}}\leq\underline{\lim}_{\epsilon\rightarrow0}\epsilon\log\mathbb{P}_{\epsilon}^{|0,\mathcal{E}}\leq\overline{\lim}_{\epsilon\rightarrow0}\epsilon\log\mathbb{P}_{\epsilon}^{|0,\mathcal{E}}\leq-U^{|0,\mathcal{E}}$$
where

$$L^{|0,\mathcal{E}}=\frac{1}{2}\inf_{h\in\xLtwo(0,+\infty;\L2):
\int_{-l}^l|\tilde{\mathbf{S}}(h)(T,x)|^2dx>4-\gamma,
\mathcal{T}\left(\tilde{\mathbf{S}}(h)\right)>t}\left\{\|h\|_{\xLtwo(0,+\infty;\L2)}^2\right\}$$
and
$$U^{|0,\mathcal{E}}=\frac{1}{2}\inf_{h\in\xLtwo(0,+\infty;\L2):\int_{-l}^l|\tilde{\mathbf{S}}(h)(T,x)|^2dx\geq
4-\gamma}\left\{\|h\|_{\xLtwo(0,+\infty;\L2)}^2 \right\}<0,$$
similarly,
$$-L^{|1,\mathcal{E}}\leq
\underline{\lim}_{\epsilon\rightarrow0}\epsilon\log\mathbb{P}_{\epsilon}^{|1,\mathcal{E}}\leq\overline{\lim}_{\epsilon\rightarrow0}\epsilon\log\mathbb{P}_{\epsilon}^{|1,\mathcal{E}}\leq-U^{|1,\mathcal{E}}$$
where

$$L^{|1,\mathcal{E}}=\frac{1}{2}\inf_{h\in\xLtwo(0,+\infty;\L2):\int_{-l}^l|\tilde{\mathbf{S}}(h)(T,x)|^2dx<4-\gamma,
\mathcal{T}\left(\mathbf{S}(h)\right)>t}\left\{\|h\|_{\xLtwo(0,+\infty;\L2)}^2
\right\}$$ and
$$U^{|1,\mathcal{E}}=\frac{1}{2}\inf_{h\in\xLtwo(0,+\infty;\L2):\int_{-l}^l|\mathbf{S}(h)(T,x)|^2dx\leq
4-\gamma}\left\{\|h\|_{\xLtwo(0,+\infty;\L2)}^2\right\}<0.$$
\end{prpstn}
\begin{proof} The result follows from the LDP
for the laws of the solutions $(u_{\epsilon}^n)_{\epsilon>0}$, the
fact that $\varphi^{-1}\left([4-\gamma,+\infty)\right)$ is a
closed set containing
$\varphi^{-1}\left([4-\gamma,+\infty)\right)\cap\mathcal{T}^{-1}\left((T,+\infty)\right)$
and that
$\varphi^{-1}\left((4-\gamma,+\infty)\right)\cap\mathcal{T}^{-1}\left((T,+\infty)\right)$
is an open set included in
$\varphi^{-1}\left([4-\gamma,+\infty)\right)\cap\mathcal{T}^{-1}\left((T,+\infty)\right)$.
\end{proof}

Remark again that, under the assumption that the null space of
$\Phi^*$ is zero, the infima are never taken over empty sets. If
$l$ is taken as $+\infty$, we obtain the same upper bounds as in
section \ref{s61}. Also, if $\Phi$ acts as the identity on a
sufficiently large linear space of $\L2$, we could implement the
previous computations of the calculus of variations on the
parametrized controls.

\end{document}